\documentclass[preprint,12pt]{elsarticle}

\usepackage[T1]{fontenc}
\usepackage{lmodern}
\usepackage[utf8]{inputenc}
\usepackage{amsmath,amssymb,amsthm,mathtools}
\usepackage[shortlabels]{enumitem}
\usepackage{url}
\usepackage[hidelinks]{hyperref}
\hypersetup{pdftitle={Exit times from time-dependent random domains: continuity, weak convergence, and exit-time profiles},pdfauthor={Tristan Guillaume}}

\biboptions{sort&compress}

\numberwithin{equation}{section}

\theoremstyle{plain}
\newtheorem{theorem}{Theorem}[section]
\newtheorem{proposition}{Proposition}[section]
\newtheorem{lemma}{Lemma}[section]
\newtheorem{corollary}{Corollary}[section]
\theoremstyle{definition}
\newtheorem{definition}{Definition}[section]
\newtheorem{example}{Example}[section]
\theoremstyle{remark}
\newtheorem{remark}{Remark}[section]

\DeclareMathOperator{\dist}{dist}
\DeclareMathOperator{\cl}{cl}
\newcommand{\NT}{\textup{(NT)}}

\newcommand{\R}{\mathbb{R}}
\newcommand{\N}{\mathbb{N}}
\newcommand{\Q}{\mathbb{Q}}
\newcommand{\E}{\mathbb{E}}
\newcommand{\PP}{\mathbb{P}}
\newcommand{\Dd}{\mathbb{D}^{d}}
\newcommand{\Done}{\mathbb{D}^{1}}
\newcommand{\cH}{\mathcal{H}}
\newcommand{\cO}{\mathcal{O}}
\newcommand{\cE}{\mathcal{E}}
\newcommand{\cN}{\mathcal{N}}
\newcommand{\Ms}{\mathbb{M}}
\newcommand{\Tcal}{T_{\mathrm{cal}}}
\newcommand{\hyp}{{}_{1}F_{1}}
\newcommand{\gbar}{\chi}

\journal{Stochastic Processes and their Applications}

\begin{document}

\begin{frontmatter}

\title{Exit times from time-dependent random domains: continuity, weak convergence, and exit-time profiles}

\author[inst1]{Tristan Guillaume\corref{cor1}}
\ead{tristan.guillaume@cyu.fr}
\cortext[cor1]{Corresponding author. CY Cergy Paris Université, Laboratoire Thema, 33 boulevard du Port, F-95011 Cergy-Pontoise Cedex, France. Tel.: +33-6-12-22-45-88.}

\affiliation[inst1]{organization={CY Cergy Paris Université, Laboratoire Thema},
            addressline={33 boulevard du Port},
            postcode={F-95011},
            city={Cergy-Pontoise Cedex},
            country={France}}

\begin{abstract}
We study exit times from time-dependent domains under joint perturbations
of the trajectory and the domain. Representing a moving domain as the
set \(\{\Phi<0\}\) for a continuous barrier \(\Phi\) on space--time, we reduce
the exit problem to a one-dimensional first-passage problem for the
scalarized path \(y(t):=\Phi(t,x(t))\) of the trajectory \(x\). Our first
main result is a deterministic continuity
theorem: the exit-time functional is continuous, under Skorokhod \(J_{1}\)
convergence of the path and locally uniform convergence of the barrier, at
every configuration satisfying an explicit non-tangency condition \NT{}. This
condition is sharp: it characterizes the continuity set of the functional. As a consequence of the continuity theorem, exit times converge in
distribution whenever paths and barriers converge jointly in distribution
and the limiting pair satisfies \NT{} almost surely; no independence or
other structural restriction between trajectory and domain is required.
Our second main result is a functional limit theorem: the profile
\(u\mapsto\tau_{u}\) of the exit times from the domains \(\{\Phi<u\}\),
viewed after right-continuous regularization as a monotone càdlàg function
on a level interval \([u_{0},u_{1}]\), converges in distribution in the
Skorokhod \(M_{1}\) topology under the same joint convergence, provided the
limiting pair satisfies \NT{} almost surely at \(u_{0}\) and \(u_{1}\); an
example shows that \(J_{1}\) convergence can fail. Concrete verification
routes for \NT{} are provided, including a non-characteristic criterion for
Itô diffusions. The framework is illustrated by two Donsker-type examples;
in the second one, the domain is random, is approximated jointly with the
trajectory, and is coupled to the state.
\end{abstract}


\end{frontmatter}

\section{Introduction}\label{sec:intro}

Exit times and boundary-crossing times are fundamental functionals of
stochastic processes: they appear in limit theorems, in pathwise stability
questions, and in probabilistic representations of boundary value
problems. In many applications one has to approximate not only the driving
path \(x\) but also the domain \(\Gamma=(\Gamma_{t})_{t\geq 0}\) itself:
interfaces move, constraints evolve in time, and the domain may be random
and coupled to the state process. This raises the following question:
when does the convergence of joint approximations \((x_{n},\Gamma^{(n)})\)
imply the convergence of the corresponding exit times?

A first difficulty is that set convergence alone is not enough. Even
strong set-theoretic convergence, such as Hausdorff or Fell~\cite{Fell1962} convergence in space--time, does not by itself control exit
or hitting times: thin approximating sets may converge to a thick limit
while remaining avoidable by a fixed trajectory. Section~\ref{sec:counterexample}
gives an explicit counterexample showing that no general continuity theory
for exit times from varying domains can be based on set convergence alone.
This is the starting point of the present paper.

Our framework is based on a barrier representation of time-dependent
domains. We write
\begin{equation}\label{eq:intro-barrier}
\Gamma_{t}=\{z\in\R^{d}:\ \Phi(t,z)<0\},\qquad t\geq 0,
\end{equation}
where \(\Phi\) is a continuous function on space--time, called a barrier.
For a càdlàg path \(x\), the scalarized path is
\begin{equation}\label{eq:intro-scalar}
y(t):=\Phi\bigl(t,x(t)\bigr),\qquad t\geq 0.
\end{equation}
The exit time from the moving domain \(\Gamma\) is then the first time at
which \(y\) reaches or exceeds the level \(0\). This reduction has two
consequences. First,
it identifies the right notion of domain approximation for
weak-convergence purposes, namely locally uniform convergence of barriers
on space--time, rather than set convergence. Second, it turns the
exit problem into a one-dimensional first-passage problem, so that the
\(J_{1}\)-stability of the paths can be combined with continuity properties
of first-passage and inverse maps for the scalarized paths.

Continuity does not hold everywhere, even in the barrier framework. The
remaining obstruction is pathwise: the scalarized limit path may approach
the boundary level without crossing it, or reach it and then graze or
stick to it, and in such cases arbitrarily small perturbations can produce
macroscopic changes in the exit time. Our first main result,
Theorem~\ref{thm:continuity}, identifies an explicit non-tangency condition
\NT{} under which the exit-time functional is continuous under joint
perturbations of the path, in the Skorokhod \(J_{1}\) topology, and of the
barrier, locally uniformly. Proposition~\ref{prop:discontinuity} shows that
\NT{} is sharp: the exit-time functional is discontinuous at every
configuration that violates \NT{}, so that \NT{} characterizes its continuity
set.

As a direct consequence of Theorem~\ref{thm:continuity}, Theorem~\ref{thm:weak} yields the weak convergence
of exit times under joint weak convergence of random paths and barriers
\((X_{n},\Phi_{n})\), provided the limiting pair satisfies \NT{} almost surely. The argument is
pathwise and requires no independence, or any other structural restriction,
between the trajectory and the domain once joint convergence is available.
It therefore applies to random environments, common-noise couplings, and
feedback-type random domains, as long as two verification tasks are
carried out: (i)~joint convergence of \((X_{n},\Phi_{n})\), and
(ii)~condition \NT{} for the limit. Section~\ref{sec:example-coupled}
carries out both tasks for a random domain that is built from the
trajectory itself.

Our second main result, Theorem~\ref{thm:profiles}, is a functional limit
theorem for exit-time profiles associated with the nested domains obtained
by varying the barrier level. For the family of domains
\begin{equation}\label{eq:intro-levels}
\Gamma^{\Phi,u}_{t}:=\{z\in\R^{d}:\ \Phi(t,z)<u\},\qquad t\geq 0,\ u\in\R,
\end{equation}
the profile \(u\mapsto\tau_{u}\) of the corresponding exit times is
nondecreasing and, after right-continuous regularization, a càdlàg function
of the level. We prove the convergence in distribution of these profiles on a
compact level interval \([u_{0},u_{1}]\), in the Skorokhod \(M_{1}\) topology,
under joint convergence of \((X_{n},\Phi_{n})\), provided the limiting pair
satisfies \NT{} almost surely at the two end levels \(u_{0}\) and \(u_{1}\). The use of \(M_{1}\) is essential: \(J_{1}\) is in general too
strong for monotone first-passage profiles, and Example~\ref{ex:J1-fails}
exhibits a situation in which \(J_{1}\) convergence fails although the
scalarized paths converge uniformly.

The role of \NT{} deserves emphasis. It is not an artifact of the proofs,
but the natural continuity-set condition for first-passage functionals. It
excludes grazing and sticking configurations, which can occur, for instance,
for degenerate diffusions whose noise vanishes in the direction normal to
the boundary, or for deterministic dynamics that slide along a moving
interface after contact.
In such cases exit-time stability may still hold, but it is necessarily
model-specific and lies outside a general continuous mapping theorem. To
make clear which time-dependent exit problems are covered by our
framework, Section~\ref{sec:verification} provides concrete verification
routes for \NT{}: a jump-overshoot criterion and an Itô criterion based on a
non-characteristic boundary condition.

\paragraph{Related literature}
The Skorokhod \(J_{1}\) and \(M_{1}\) topologies were introduced by
Skorokhod~\cite{Skorokhod1956}; see Billingsley~\cite{Billingsley1999},
Jacod and Shiryaev~\cite{JacodShiryaev2003} and Whitt~\cite{Whitt2002} for
modern treatments. The continuity of first-passage and inverse maps under
Skorokhod convergence goes back to Whitt~\cite{Whitt1971,Whitt1980,Whitt2002}
and Vervaat~\cite{Vervaat1972}, and is a standard tool in weak-convergence
theory; the \(M_{1}\) topology is the natural setting for monotone and
inverse processes, and its role for first-passage-time processes is well
documented~\cite{PuhalskiiWhitt1997,Whitt2002}. The present work combines
this pathwise viewpoint with a barrier (level-set) representation of
time-dependent domains and a locally uniform topology on barriers, which
yields a joint robustness framework for varying, possibly random, domains.
The barrier formulation is in line with level-set methods and
viscosity-solution numerics~\cite{BarlesSouganidis1991,CrandallIshiiLions1992,GibouFedkiwOsher2018,OsherSethian1988},
in which moving interfaces are typically approximated through locally
uniformly convergent level-set functions. Time-dependent domains also
appear prominently in the literature on reflected processes and on
probabilistic representations of parabolic boundary value problems; see,
for example,~\cite{NystromOnskog2010} for the Skorokhod problem in
time-dependent domains. In the one-dimensional continuous setting,
change-of-variable formulas with local time on curves provide another
classical route to boundary crossing and moving interfaces~\cite{Peskir2005}.

The framework also connects with several recent lines of work. In
discrete-time optimal stopping, Soner and
Tissot-Daguette~\cite{SonerTissotDaguette2025} study relaxations of the
hitting times of stopping boundaries and the continuity, with respect to
the boundary, of the expected reward of the associated stopping rule; the
continuity of optimal stopping surfaces for jump-diffusions is studied by
Cai, De Angelis and Palczewski~\cite{CaiDeAngelisPalczewski2023}. These
works are motivated by optimal stopping, where the relevant boundaries are
determined by value functions, and they concern expected rewards or optimal
boundaries, whereas the present paper treats prescribed, possibly random
and state-coupled, domains and the convergence of the exit times
themselves. In a different direction, Liang and
Borovkov~\cite{LiangBorovkov2025} establish the differentiability of
boundary-crossing probabilities of diffusions with respect to the
boundary, a smooth-dependence counterpart of the weak-convergence results
obtained here.

\paragraph{Organization of the paper}
Section~\ref{sec:prelim} introduces the barrier representation, the
relevant topologies and the notion of regular level, and ends with the
counterexample showing why set convergence alone is insufficient.
Section~\ref{sec:continuity} contains the deterministic continuity theorem
for exit times, its sharpness, and its weak-convergence consequence,
together with their extensions to exit positions and stopped paths.
Section~\ref{sec:profiles} proves the functional limit theorem for
exit-time profiles in the \(M_{1}\) topology,
together with its finite-dimensional and open-interval versions,
and explains why \(M_{1}\) rather
than \(J_{1}\) is needed. Section~\ref{sec:verification} provides
verification routes for
\NT{} and two worked Donsker-type examples: a deterministic moving boundary
(Section~\ref{sec:example-donsker}) and a random tracking domain that is
approximated jointly with the trajectory and coupled to it
(Section~\ref{sec:example-coupled}). \hyperref[app:M1]{Appendix~A} contains the
proof of the monotone \(M_{1}\) criterion used for the profiles.

\section{Preliminaries: barrier domains, topologies, and first-passage maps}\label{sec:prelim}

\subsection{Path spaces and topologies}\label{sec:topologies}

Fix \(d\geq 1\). For \(k\geq 1\), let \(\mathbb{D}^{k}:=D([0,\infty),\R^{k})\) be
the space of càdlàg functions \(x:[0,\infty)\to\R^{k}\); we mainly work with
\(\Dd\) and with the space \(\Done\) of real càdlàg functions.

Throughout, \(\N:=\{1,2,\ldots\}\), \(\|\cdot\|\) denotes the Euclidean norm on \(\R^{m}\), \(m\geq 1\),
and \(|\cdot|\) the absolute value on \(\R\); \(\inf\varnothing:=+\infty\), and
\([0,\infty]\) carries its usual compact topology. Lowercase letters
\(x,y\) denote deterministic paths, while uppercase letters \(X,Y\) denote
random elements of the same spaces (processes). For a càdlàg function \(f\)
with values in \(\R^{m}\), we write \(f(t-)\) for its left limit at \(t>0\),
set \(f(0-):=f(0)\), and write \(\Delta f(t):=f(t)-f(t-)\) for its jump at
\(t\). For topological spaces \(A\) and \(B\), \(C(A,B)\) is the space of
continuous maps from \(A\) to \(B\); \(C^{1}([0,\infty))\) is the space of
continuously differentiable functions on \([0,\infty)\), and \(C^{1,2}\)
denotes the class of functions of \((t,z)\) that are continuously
differentiable once in \(t\) and twice in \(z\).

We equip each \(\mathbb{D}^{k}\) with the
Skorokhod \(J_{1}\) topology for paths indexed by the half-line (sometimes
called the local \(J_{1}\) topology); see~\cite[Section~16]{Billingsley1999},~\cite[Section~VI.1]{JacodShiryaev2003} and~\cite{Whitt2002}. Let \(\Lambda\)
be the set of continuous, strictly increasing functions
\(\lambda:[0,\infty)\to[0,\infty)\) with \(\lambda(0)=0\) and
\(\lambda(t)\to\infty\) as \(t\to\infty\). We shall repeatedly use the
following characterization~\cite[Section~VI.1]{JacodShiryaev2003}:
\(x_{n}\to x\) in \(\Dd\) if and only if there exist \(\lambda_{n}\in\Lambda\)
such that
\begin{equation}\label{eq:J1-char}
\sup_{t\geq 0}|\lambda_{n}(t)-t|\to 0,\qquad
\sup_{0\leq t\leq N}\|x_{n}(\lambda_{n}(t))-x(t)\|\to 0\quad(N\in\N).
\end{equation}
Equivalently, the restrictions of \(x_{n}\) to \([0,T]\) converge to the
restriction of \(x\) in the \(J_{1}\) topology of \(D([0,T],\R^{d})\) for every \(T>0\) at
which \(x\) is continuous~\cite[Section~16]{Billingsley1999}. With this
topology, each \(\mathbb{D}^{k}\) is a Polish space whose Borel \(\sigma\)-field is generated
by the evaluation maps \(x\mapsto x(t)\), \(t\geq 0\)~\cite[Section~16]{Billingsley1999}. Two elementary consequences of
\eqref{eq:J1-char} will be used: if \(x_{n}\to x\) in \(\Dd\), then
\(\sup_{n}\sup_{t\leq T}\|x_{n}(t)\|<\infty\) for every \(T>0\), and
\(x_{n}(t)\to x(t)\) at every continuity point \(t\) of \(x\).

Let \(E:=[0,\infty)\times\R^{d}\subset\R^{1+d}\) denote space--time; for
\(p,q\in E\), we write \(|p-q|:=\|p-q\|\) for their Euclidean distance in \(\R^{1+d}\). Let
\(\mathbb{B}_{R}:=\{z\in\R^{d}:\|z\|\leq R\}\) be the closed ball of radius
\(R\) centered at the origin. Barriers are elements of
\[
\cH:=C(E,\R),
\]
equipped with the topology of uniform convergence on compact subsets of
\(E\):
\begin{equation}\label{eq:H-conv}
\Phi_{n}\to\Phi\ \text{in}\ \cH
\quad\Longleftrightarrow\quad
\sup_{[0,T]\times\mathbb{B}_{R}}|\Phi_{n}-\Phi|\to 0\quad(T,R>0).
\end{equation}
This topology is induced by the complete and separable metric
\[
d_{\cH}(\Phi,\Psi):=\sum_{k\geq 1}2^{-k}\min\Bigl\{1,\sup_{[0,k]\times\mathbb{B}_{k}}|\Phi-\Psi|\Bigr\},
\]
so that \(\cH\) and \(\Dd\times\cH\) are Polish spaces.

\subsection{Time-dependent domains, barrier representation, and scalarization}\label{sec:barrier}

A time-dependent open domain is a family \(\Gamma=(\Gamma_{t})_{t\geq 0}\)
of open sets \(\Gamma_{t}\subset\R^{d}\). Its admissible and forbidden
space--time sets are
\begin{equation}\label{eq:OF}
\cO(\Gamma):=\{(t,z)\in E:\ z\in\Gamma_{t}\},\qquad
F(\Gamma):=E\setminus\cO(\Gamma).
\end{equation}
The exit time of a path \(x\in\Dd\) from \(\Gamma\) is
\begin{equation}\label{eq:exit-Gamma}
\tau(x,\Gamma):=\inf\{t\geq 0:\ x(t)\notin\Gamma_{t}\}
=\inf\{t\geq 0:\ (t,x(t))\in F(\Gamma)\}.
\end{equation}
For \(\Phi\in\cH\) and \(u\in\R\), the barrier domain
\(\Gamma^{\Phi,u}\) and the associated exit time are
\begin{equation}\label{eq:barrier-domain}
\begin{aligned}
\Gamma^{\Phi,u}_{t}&:=\{z\in\R^{d}:\ \Phi(t,z)<u\},\qquad t\geq 0,\\
\tau_{u}(x,\Phi)&:=\inf\{t\geq 0:\ \Phi(t,x(t))\geq u\},
\end{aligned}
\end{equation}
so that \(\tau_{u}(x,\Phi)=\tau(x,\Gamma^{\Phi,u})\). We write
\(\Gamma^{\Phi}:=\Gamma^{\Phi,0}\) and \(\tau(x,\Phi):=\tau_{0}(x,\Phi)\).

\begin{remark}[Only continuity of \(\Phi\) is used]\label{rem:continuity-only}
No differentiability of \(\Phi\) is required in
Sections~\ref{sec:prelim}--\ref{sec:profiles}. In particular, the framework
allows nonsmooth moving domains, with corners or edges, and barriers
constructed by maximum and minimum operations; for instance, a moving
convex polyhedral domain can be encoded by a finite maximum of functions
\((t,z)\mapsto\alpha(t)\cdot z-\beta(t)\) with continuous coefficients
\(\alpha\) and \(\beta\), and a general polyhedral domain by finite max--min
combinations of such functions.
\end{remark}

The next lemma characterizes the domains that admit a continuous barrier
and provides a canonical choice. For \(A\subset E\), \(\cl A\) and
\(\partial A\) denote the closure and the boundary of \(A\) relative to \(E\),
and \(\dist(p,A):=\inf_{q\in A}|p-q|\) for \(p\in E\).

\begin{lemma}[Universality of the barrier representation]\label{lem:universality}
Let \(\Gamma\) be a time-dependent open domain. The following are
equivalent:
\begin{enumerate}[(i)]
\item \(\cO(\Gamma)\) is open in \(E\) (equivalently, \(F(\Gamma)\) is closed
in \(E\));
\item there exists \(\Phi\in\cH\) such that \(\Gamma_{t}=\Gamma^{\Phi}_{t}\) for
all \(t\geq 0\).
\end{enumerate}
Moreover, if (i) holds and \(\cO(\Gamma)\) and \(F(\Gamma)\) are nonempty,
then the signed distance
\begin{equation}\label{eq:signed-distance}
\Phi_{\Gamma}(p):=\dist\bigl(p,\cO(\Gamma)\bigr)-\dist\bigl(p,F(\Gamma)\bigr),\qquad p\in E,
\end{equation}
is \(1\)-Lipschitz and satisfies
\[
\{\Phi_{\Gamma}<0\}=\cO(\Gamma),\qquad
\{\Phi_{\Gamma}=0\}=\partial\cO(\Gamma),\qquad
\{\Phi_{\Gamma}>0\}=E\setminus\cl\cO(\Gamma).
\]
In particular, \(\Phi_{\Gamma}\in\cH\) represents \(\Gamma\) in the sense
of (ii). (If \(F(\Gamma)=\varnothing\), respectively
\(\cO(\Gamma)=\varnothing\), the constant barrier \(-1\), respectively
\(1\), represents \(\Gamma\).)
\end{lemma}

\begin{proof}
If \(\Gamma=\Gamma^{\Phi}\) with \(\Phi\in\cH\), then
\(\cO(\Gamma)=\Phi^{-1}((-\infty,0))\) is open in \(E\). Conversely, assume
that \(\cO:=\cO(\Gamma)\) is open, so that \(F:=F(\Gamma)\) is closed, and
that both sets are nonempty. If \(p\in\cO\), then \(\dist(p,\cO)=0\), while
\(\dist(p,F)>0\) because \(F\) is closed and \(p\notin F\); hence
\(\Phi_{\Gamma}(p)<0\). If \(p\in F\), then \(\Phi_{\Gamma}(p)=\dist(p,\cO)\geq 0\),
with equality if and only if \(p\in\cl\cO\). This gives the three identities,
since \(\partial\cO=\cl\cO\setminus\cO=F\cap\cl\cO\).

For the Lipschitz property, let \(p,q\in E\). If \(p,q\in\cO\) or
\(p,q\in F\), then \(|\Phi_{\Gamma}(p)-\Phi_{\Gamma}(q)|\leq|p-q|\), because
distance functions are \(1\)-Lipschitz. In the remaining case, one of the
two points lies in \(\cO\) and the other in \(F\); by symmetry, we may assume
that \(p\in\cO\) and \(q\in F\). Since \(E\) is convex, the points \(p+s(q-p)\), \(s\in[0,1]\),
lie in \(E\). Let \(s^{*}:=\sup\{{s\in[0,1]:\ p+s(q-p)\in\cO}\}\) and
\(w:=p+s^{*}(q-p)\). Then \(w\in\cl\cO\). Moreover, \(w\notin\cO\):
otherwise \(s^{*}<1\), because \(q\notin\cO\), and since \(\cO\) is open in
\(E\), we would have \(p+s(q-p)\in\cO\) for all \(s\) slightly larger than
\(s^{*}\), which contradicts the definition of \(s^{*}\). Hence \(w\in F\),
and
\[
0<\Phi_{\Gamma}(q)-\Phi_{\Gamma}(p)=\dist(q,\cO)+\dist(p,F)\leq|q-w|+|w-p|=|q-p|,
\]
because \(w\) lies on the segment between \(p\) and \(q\).
\end{proof}

For \(f\in\Done\) and \(u\in\R\), we define the first-passage times
\begin{equation}\label{eq:first-passage}
T_{f}(u):=\inf\{t\geq 0:\ f(t)\geq u\},\qquad
T^{+}_{f}(u):=\inf\{t\geq 0:\ f(t)>u\}.
\end{equation}
Both are nondecreasing functions of \(u\), and \(T_{f}\leq T^{+}_{f}\).

\begin{corollary}[Scalarization]\label{cor:scalarization}
Let \(\Phi\in\cH\) and \(x\in\Dd\), and define the scalarized path
\(y(t):=\Phi(t,x(t))\), \(t\geq 0\). Then \(y\in\Done\) and
\(\tau_{u}(x,\Phi)=T_{y}(u)\) for every \(u\in\R\). In particular,
\(\tau(x,\Gamma^{\Phi})=\tau(x,\Phi)=T_{y}(0)\).
\end{corollary}

\begin{proof}
Since \(t\mapsto(t,x(t))\) is càdlàg and \(\Phi\) is continuous, \(y\) is
càdlàg. The identity \(\tau_{u}(x,\Phi)=T_{y}(u)\) is the definition
\eqref{eq:barrier-domain} rewritten in terms of \(y\).
\end{proof}

\begin{lemma}[Continuity of the scalarization map]\label{lem:scalarization}
Let \(x_{n}\to x\) in \(\Dd\) and \(\Phi_{n}\to\Phi\) in \(\cH\), and let
\(\lambda_{n}\in\Lambda\) satisfy \eqref{eq:J1-char}. Define
\(y_{n}(t):=\Phi_{n}(t,x_{n}(t))\) and \(y(t):=\Phi(t,x(t))\). Then
\begin{equation}\label{eq:scalar-conv}
\sup_{0\leq t\leq N}|y_{n}(\lambda_{n}(t))-y(t)|\to 0\qquad\text{for every } N\in\N .
\end{equation}
In particular, \(y_{n}\to y\) in \(\Done\), and the scalarization map
\(S:\Dd\times\cH\to\Done\), \(S(x,\Phi):=\Phi(\cdot,x(\cdot))\), is
continuous.
\end{lemma}

\begin{proof}
Fix \(N\in\N\) and let \(R:=1+\sup_{t\leq N}\|x(t)\|\), which is finite
because càdlàg functions are bounded on compact intervals. The set
\(K:=[0,N+1]\times\mathbb{B}_{R}\) is compact, and by \eqref{eq:J1-char} there
is \(n_{0}\) such that \((\lambda_{n}(t),x_{n}(\lambda_{n}(t)))\in K\) for all
\(n\geq n_{0}\) and \(t\in[0,N]\); clearly \((t,x(t))\in K\) as well. Let
\(\omega_{K}(\delta):=\sup\{|\Phi(p)-\Phi(q)|:\ p,q\in K,\ |p-q|\leq\delta\}\)
be the modulus of continuity of \(\Phi\) on \(K\), so that
\(\omega_{K}(\delta)\to 0\) as \(\delta\downarrow 0\), and let
\[
\delta_{n}:=\sup_{s\geq 0}|\lambda_{n}(s)-s|+\sup_{0\leq s\leq N}\|x_{n}(\lambda_{n}(s))-x(s)\| .
\]
For \(n\geq n_{0}\) and \(t\in[0,N]\), writing
\(p_{n}(t):=(\lambda_{n}(t),x_{n}(\lambda_{n}(t)))\in K\), we obtain
\[
\begin{aligned}
|y_{n}(\lambda_{n}(t))-y(t)|
&\leq\bigl|\Phi_{n}(p_{n}(t))-\Phi(p_{n}(t))\bigr|+\bigl|\Phi(p_{n}(t))-\Phi\bigl(t,x(t)\bigr)\bigr|\\
&\leq\sup_{K}|\Phi_{n}-\Phi|+\omega_{K}(\delta_{n}).
\end{aligned}
\]
Both terms tend to \(0\) by \eqref{eq:H-conv} and \eqref{eq:J1-char}, which
proves \eqref{eq:scalar-conv}. Since the same \(\lambda_{n}\) serve for
every \(N\), the characterization \eqref{eq:J1-char} (with \(d=1\)) yields
\(y_{n}\to y\) in \(\Done\).
\end{proof}

\begin{remark}[Scope of the barrier framework]\label{rem:scope}
By Lemma~\ref{lem:universality}, the domains that admit a continuous
barrier are exactly those whose admissible space--time set \(\cO(\Gamma)\)
is open. Openness of each slice \(\Gamma_{t}\) is necessary for this, but not
sufficient. For instance, if \(d=1\) and \(\Gamma_{t}=(-\infty,g(t))\) for a
function \(g:[0,\infty)\to\R\), then \(\cO(\Gamma)\) is open if and only if
\(g\) is lower semicontinuous. For a càdlàg \(g\), this means that all jumps
of \(g\) are downward. Thus, in this one-sided model, domains that shrink
discontinuously in time are covered, although the natural barrier
\((t,z)\mapsto z-g(t)\) is then discontinuous and has to be replaced, for
instance, by \(\Phi_{\Gamma}\) or by an augmented formulation (see
Remark~\ref{rem:drawdown}), whereas domains that expand discontinuously in
time lie outside the present framework. For general families of open sets,
the criterion is the openness of \(\cO(\Gamma)\); the distinction between
abrupt shrinking and abrupt expansion is only a reading of it in one-sided
models such as the one above. Domains for which \(\cO(\Gamma)\) is not open
would require a different space and topology for the domain variable, for
instance a space of barriers that are only càdlàg in time.
\end{remark}

\subsection{Regular levels}\label{sec:regular}

The continuity of first-passage maps is governed by the following standard
notion; cf.~\cite{Whitt1980,Whitt2002}.

\begin{definition}[Regular level]\label{def:regular}
Let \(f\in\Done\) and \(u\in\R\) with \(T_{f}(u)<\infty\). The level \(u\) is
\emph{regular} for \(f\) if
\begin{align}
&\sup_{0\leq s\leq T}f(s)<u\quad\text{for every } T\in[0,T_{f}(u)),\tag{R$-$}\label{eq:Rminus}\\
&\sup_{T_{f}(u)\leq s\leq T_{f}(u)+\varepsilon}f(s)>u\quad\text{for every } \varepsilon>0 .\tag{R$+$}\label{eq:Rplus}
\end{align}
\end{definition}

Condition \eqref{eq:Rplus} says that \(f\) strictly overshoots \(u\)
arbitrarily soon after first reaching it. It fails exactly when there is
\(\varepsilon>0\) such that \(f(s)\leq u\) for all
\(s\in[T_{f}(u),T_{f}(u)+\varepsilon]\), that is, when \(f\) does not exceed \(u\) on some time
interval of positive length starting at \(T_{f}(u)\), as for paths that
stick to \(u\). Returning below \(u\) is compatible with \eqref{eq:Rplus},
provided strict exceedances of \(u\) still occur arbitrarily soon after
\(T_{f}(u)\), as they do for Brownian paths at a fixed level, almost
surely. Condition \eqref{eq:Rminus}
says that \(f\) stays below \(u\) uniformly on every compact subinterval of
\([0,T_{f}(u))\). Since \(f(s)<u\) for all \(s<T_{f}(u)\), \eqref{eq:Rminus}
can only fail if \(f\) approaches \(u\) through a left limit, i.e.\ if
\(f(t-)=u\) for some \(t\in(0,T_{f}(u))\) (see Remark~\ref{rem:NTminus}); in
particular, it holds automatically when \(f\) is continuous.

Regularity has a simple interpretation in terms of the first-passage
profile \(v\mapsto T_{f}(v)\).

\begin{lemma}[Regular levels and continuity of the first-passage profile]\label{lem:regular-continuity}
Let \(f\in\Done\) and \(u\in\R\).
\begin{enumerate}[(a)]
\item \(T^{+}_{f}(u)=\lim_{v\downarrow u}T_{f}(v)\). Consequently,
\(T^{+}_{f}\) is the right-continuous version of the nondecreasing function
\(T_{f}\).
\item If \(T_{f}(u)<\infty\), then \eqref{eq:Rplus} holds if and only if
\(T^{+}_{f}(u)=T_{f}(u)\), and \eqref{eq:Rminus} holds if and only if
\(\lim_{v\uparrow u}T_{f}(v)=T_{f}(u)\).
\end{enumerate}
In particular, if \(T_{f}(u)<\infty\), the level \(u\) is regular for \(f\)
if and only if \(v\mapsto T_{f}(v)\) is continuous at \(v=u\), and then
\(T^{+}_{f}(u)=T_{f}(u)\).
\end{lemma}

\begin{proof}
(a) For \(v>u\), every \(t\) with \(f(t)\geq v\) satisfies \(f(t)>u\), so
\(T_{f}(v)\geq T^{+}_{f}(u)\); the limit \(\lim_{v\downarrow u}T_{f}(v)\)
exists by monotonicity and is therefore at least \(T^{+}_{f}(u)\).
Conversely, if \(f(t)>u\), then \(T_{f}(v)\leq t\) for all
\(v\in(u,f(t)]\), so that \(\lim_{v\downarrow u}T_{f}(v)\leq t\); taking the
infimum over all such \(t\) gives \(\lim_{v\downarrow u}T_{f}(v)\leq T^{+}_{f}(u)\).
The right limit of a nondecreasing function is right-continuous.

(b) Let \(\tau:=T_{f}(u)<\infty\). Since \(f(s)>u\) implies \(f(s)\geq u\), we
have \(T^{+}_{f}(u)\geq\tau\). If \eqref{eq:Rplus} holds, then for every
\(\varepsilon>0\) there is \(s\in[\tau,\tau+\varepsilon]\) with \(f(s)>u\),
so \(T^{+}_{f}(u)\leq\tau+\varepsilon\); hence \(T^{+}_{f}(u)=\tau\).
Conversely, if \(T^{+}_{f}(u)=\tau\), then for every \(\varepsilon>0\) there
is \(s<\tau+\varepsilon\) with \(f(s)>u\), and necessarily \(s\geq\tau\);
hence \eqref{eq:Rplus} holds.

Suppose next that \eqref{eq:Rminus} holds and fix \(T\in[0,\tau)\). Then
\(m_{T}:=\sup_{[0,T]}f<u\), and for \(v\in(m_{T},u)\) we have \(f(t)<v\) for
all \(t\leq T\), so that \(T_{f}(v)\geq T\). Hence
\(\lim_{v\uparrow u}T_{f}(v)\geq T\) for every \(T<\tau\), and since
\(T_{f}(v)\leq\tau\) for \(v\leq u\), the left limit equals \(\tau\).
Conversely, if \eqref{eq:Rminus} fails, there is \(T<\tau\) with
\(\sup_{[0,T]}f\geq u\). Then, for every \(v<u\), there is \(t\leq T\) with
\(f(t)>v\), so that \(T_{f}(v)\leq T\), and
\(\lim_{v\uparrow u}T_{f}(v)\leq T<\tau\).
\end{proof}

\subsection{Why set convergence alone does not control exit times: a counterexample}\label{sec:counterexample}

Random time-dependent domains are naturally described through random
closed sets, for instance their forbidden sets \(F(\Gamma)\), and random
closed sets are commonly studied through the Fell (hit-or-miss) topology~\cite{Fell1962,Matheron1975,Molchanov2005}. While Fell- and Hausdorff-type
convergences are powerful tools for geometric functionals, they do not in
general control hitting or exit times: a fixed trajectory may avoid every
set of a sequence of thin sets that converges to a thick limit. The
following example makes this precise with time-independent domains in
\(\R^{2}\). For a closed set \(A\subset\R^{2}\) and a continuous path
\(\gamma\), the hitting time \(H_{A}(\gamma):=\inf\{t\geq 0:\ \gamma(t)\in A\}\)
is the exit time of \(\gamma\) from the open domain \(\R^{2}\setminus A\).
For nonempty compact sets \(A,A'\subset\R^{2}\), the Hausdorff distance is
\[
d_{\mathrm{Haus}}(A,A'):=\max\Bigl\{\sup_{a\in A}\inf_{a'\in A'}\|a-a'\|,\ \sup_{a'\in A'}\inf_{a\in A}\|a-a'\|\Bigr\}.
\]

Let \(K:=\{z\in\R^{2}:\ \|z\|\leq 1\}\) be the closed unit disk and, for
\(n\geq 1\), let
\[
K_{n}:=\bigcup_{k=0}^{n}\Bigl\{z\in\R^{2}:\ \|z\|=\frac{k}{n}\Bigr\}
\]
be the union of the concentric circles of radii \(k/n\), \(0\leq k\leq n\).
Then \(K_{n}\subset K\). Conversely, if \(z\in K\), choose
\(k\in\{0,\ldots,n\}\) with \(\bigl|\|z\|-k/n\bigr|\leq 1/n\) and let \(z'\) be
the point of norm \(k/n\) on the ray through \(z\) (with \(z':=0\) if
\(z=0\)); then \(z'\in K_{n}\) and \(\|z-z'\|=\bigl|\|z\|-k/n\bigr|\leq 1/n\).
Hence \(d_{\mathrm{Haus}}(K_{n},K)\leq 1/n\to 0\), and in particular \(K_{n}\to K\) in the Fell
topology.

Now fix an irrational number \(r\in(0,1)\), and let \(\gamma\) be the path
that runs along the circle of radius \(r\), namely
\(\gamma(t):=(r\cos t,\,r\sin t)\) for \(t\geq 0\). Since \(r\notin\{k/n:\ 0\leq k\leq n\}\) for every \(n\), the
path never meets \(K_{n}\), whereas \(\gamma(0)\in K\). Therefore
\[
H_{K_{n}}(\gamma)=\infty\quad\text{for all } n,\qquad H_{K}(\gamma)=0 .
\]
Thus a purely set-theoretic mode of domain convergence, even Hausdorff
convergence of compact sets, cannot yield a general continuity theory for
hitting or exit times. In contrast, representing domains by continuous
barriers and imposing locally uniform convergence in \(\cH\) enforces the
required geometric ``thickness''. Under barrier convergence, the remaining
obstruction is pathwise (grazing or sticking versus non-tangential
crossing), and this is precisely what condition \NT{} in
Section~\ref{sec:continuity} isolates.

\section{Continuity and weak convergence of exit times}\label{sec:continuity}

We establish the continuity of the exit-time functional under joint
perturbations of the path \(x\), in the Skorokhod \(J_{1}\) topology, and of
the barrier \(\Phi\), in \(\cH\). We then deduce the weak convergence of exit
times.
Both results are complemented by their counterparts for the exit position
and for the path stopped at exit (Corollaries~\ref{cor:exit-position}
and~\ref{cor:joint}).
Throughout this section, \((x,\Phi)\in\Dd\times\cH\) is a
\emph{configuration}, \(y(t):=\Phi(t,x(t))\), and
\[
\tau(x,\Phi)=T_{y}(0)=\inf\{t\geq 0:\ y(t)\geq 0\}\in[0,\infty]
\]
is the exit time from the barrier domain \(\Gamma^{\Phi}\)
(Corollary~\ref{cor:scalarization}). We state our results at level \(0\);
the case of a general level \(u\) is obtained by replacing \(\Phi\) with
\(\Phi-u\) (Definition~\ref{def:NT} and Remark~\ref{rem:levels}).

\subsection{A non-tangency condition}\label{sec:NT}

\begin{definition}[Condition \NT{}: non-tangential crossing]\label{def:NT}
The configuration \((x,\Phi)\) satisfies condition \NT{} if
\begin{equation*}\tag{NT$-$}\label{eq:NTminus}
\sup_{0\leq s\leq T}y(s)<0\qquad\text{for every finite } T\in[0,\tau(x,\Phi)),
\end{equation*}
and if, in the case \(\tau(x,\Phi)<\infty\),
\begin{equation*}\tag{NT$+$}\label{eq:NTplus}
\sup_{\tau(x,\Phi)\leq s\leq\tau(x,\Phi)+\varepsilon}y(s)>0\qquad\text{for every } \varepsilon>0 .
\end{equation*}
More generally, for \(u\in\R\), the configuration \((x,\Phi)\) satisfies \NT{}
(respectively \eqref{eq:NTminus}, \eqref{eq:NTplus}) \emph{at level} \(u\) if
\((x,\Phi-u)\) satisfies \NT{} (respectively \eqref{eq:NTminus},
\eqref{eq:NTplus}).
\end{definition}

In words, the scalarized path stays uniformly below the boundary level on
every compact time interval before the exit time and, if exit occurs, it
strictly overshoots the boundary level arbitrarily soon after first
reaching it. If \(\tau(x,\Phi)<\infty\), \NT{} says exactly that the level
\(0\) is regular for \(y\) in the sense of Definition~\ref{def:regular}. If
\(\tau(x,\Phi)=\infty\), \NT{} reduces to \eqref{eq:NTminus}, i.e.\
\(\sup_{[0,T]}y<0\) for every \(T\geq 0\).

\begin{remark}[On \eqref{eq:NTminus}]\label{rem:NTminus}
For every càdlàg function \(f\) and every \(T\geq 0\),
\[
\sup_{0\leq s\leq T}f(s)=\max_{0\leq t\leq T}\max\{f(t),f(t-)\}.
\]
The supremum is thus attained, either as a value or as a left limit, but it need
not be attained as a value. Since \(y(t)<0\) for every \(t<\tau(x,\Phi)\), it
follows that \eqref{eq:NTminus} is equivalent to
\begin{equation}\label{eq:NTminus-equiv}
y(t-)<0\qquad\text{for every } t\in\bigl(0,\tau(x,\Phi)\bigr).
\end{equation}
In particular, \eqref{eq:NTminus} holds automatically if \(y\) has no
downward jumps on \((0,\tau(x,\Phi))\), for instance if \(y\) is continuous,
since then \(y(t-)\leq y(t)<0\). It is not automatic in general, even when
\(\tau(x,\Phi)<\infty\). For example, let \(d=1\), \(\Phi(t,z):=z\), and
\(x(t):=t-\lfloor t\rfloor-1\) for \(t<2\), \(x(t):=1\) for \(t\geq 2\).
Then \(y=x<0\) on \([0,2)\) and \(\tau(x,\Phi)=2\), but
\(\sup_{0\leq s\leq T}y(s)=0\) for every \(T\in[1,2)\): the path approaches
the boundary through the left limit \(y(1-)=0\) without reaching it.
Proposition~\ref{prop:discontinuity} shows that the exit-time functional is
discontinuous at such configurations.
\end{remark}

\begin{remark}[The case \(\tau(x,\Phi)=\infty\)]\label{rem:no-exit}
The regime \(\tau(x,\Phi)=\infty\) is not pathological and may be generic,
for instance for expanding domains, for dynamics that are transient
relative to the boundary, or in trapping regions. In the barrier
formulation, \(\tau(x,\Phi)=\infty\) means that \(y(t)<0\) for all \(t\).
Continuity at such configurations requires more than this pointwise strict
negativity, namely the uniform interior buffer \eqref{eq:NTminus} on every
compact time interval. If it fails, the trajectory comes arbitrarily close
to the boundary on some finite horizon, and arbitrarily small perturbations
of the barrier create a finite exit time before that
horizon; see Proposition~\ref{prop:discontinuity}.
\end{remark}

\begin{remark}[Dependence of \NT{} on the representation, and the signed-distance barrier]\label{rem:representation}
Condition \NT{} depends on the barrier \(\Phi\), and not only on the domain
\(\Gamma^{\Phi}\). It is invariant under \(\Phi\mapsto h\circ\Phi\) for every
continuous, strictly increasing \(h:\R\to\R\) with \(h(0)=0\), but not under
arbitrary changes of representation. For instance, \(\min\{\Phi,0\}\)
represents the same domain as \(\Phi\), but \eqref{eq:NTplus} fails for it
whenever \(\tau(x,\Phi)<\infty\), since the corresponding scalarized path
never becomes positive. The
signed-distance barrier \(\Phi_{\Gamma}\) of Lemma~\ref{lem:universality}
avoids this degeneracy: its zero set is \(\partial\cO(\Gamma)\), and it is
positive exactly outside \(\cl\cO(\Gamma)\). For \(\Phi=\Phi_{\Gamma}\),
condition \eqref{eq:NTminus} says that the lifted path \((t,x(t))\) stays at
a positive distance from \(F(\Gamma)\) on every compact interval
\([0,T]\subset[0,\tau(x,\Gamma))\), and \eqref{eq:NTplus} says that it enters
the exterior \(E\setminus\cl\cO(\Gamma)\) at times \(t\geq\tau(x,\Gamma)\)
arbitrarily close to \(\tau(x,\Gamma)\). Both are intrinsic geometric properties of the pair
\((x,\Gamma)\). In practice, one should therefore work with barriers whose
zero set is not artificially thickened. Note also that \(\Phi_{\Gamma}\) is
Lipschitz but in general not \(C^{1,2}\), so that the Itô-type criteria of
Section~\ref{sec:routeB} do not apply to it directly; for those criteria
one uses a smoother barrier representing the same domain, which is
available, for instance, for sufficiently regular boundaries. The results
of Sections~\ref{sec:prelim}--\ref{sec:profiles} only require continuity of
the barrier.
\end{remark}

\subsection{Deterministic continuity of the exit-time functional}\label{sec:continuity-thm}

\begin{theorem}[Deterministic continuity]\label{thm:continuity}
Let \(x_{n}\to x\) in \(\Dd\) and \(\Phi_{n}\to\Phi\) in \(\cH\). If \((x,\Phi)\)
satisfies \NT{}, then
\[
\tau(x_{n},\Phi_{n})\to\tau(x,\Phi)\qquad\text{in } [0,\infty].
\]
\end{theorem}

\begin{proof}
Let \(\lambda_{n}\in\Lambda\) satisfy \eqref{eq:J1-char}, write
\(y_{n}(t):=\Phi_{n}(t,x_{n}(t))\) and \(y(t):=\Phi(t,x(t))\), and let
\[
\eta_{n}(N):=\sup_{0\leq t\leq N}|y_{n}(\lambda_{n}(t))-y(t)|,\qquad N\in\N,
\]
so that \(\eta_{n}(N)\to 0\) for every \(N\) by
Lemma~\ref{lem:scalarization}. Write \(\tau:=\tau(x,\Phi)\),
\(\tau^{(n)}:=\tau(x_{n},\Phi_{n})\), and
\[
\sigma_{n}:=\inf\{t\geq 0:\ y_{n}(\lambda_{n}(t))\geq 0\}.
\]
Since \(\lambda_{n}\) is an increasing bijection of \([0,\infty)\), we have
\(\tau^{(n)}=\lambda_{n}(\sigma_{n})\), with the convention
\(\lambda_{n}(\infty):=\infty\). Because
\(\sup_{t}|\lambda_{n}(t)-t|\to 0\), it suffices to show that
\(\sigma_{n}\to\tau\) in \([0,\infty]\).

\emph{Lower bound.} Let \(T\in[0,\tau)\) be finite; if \(\tau=0\) there is
nothing to prove. By \eqref{eq:NTminus}, \(m:=\sup_{[0,T]}y<0\). Let
\(N\geq T\) be an integer. For all \(n\) such that \(\eta_{n}(N)<|m|\),
we have \(y_{n}(\lambda_{n}(t))\leq y(t)+\eta_{n}(N)<0\) for every
\(t\in[0,T]\), hence \(\sigma_{n}\geq T\). Therefore
\(\liminf_{n}\sigma_{n}\geq T\) for every finite \(T<\tau\), that is,
\(\liminf_{n}\sigma_{n}\geq\tau\). If \(\tau=\infty\), this already shows that
\(\sigma_{n}\to\infty\).

\emph{Upper bound.} Assume that \(\tau<\infty\) and let \(\varepsilon>0\).
By \eqref{eq:NTplus} there is \(t\in[\tau,\tau+\varepsilon]\) with
\(y(t)>0\). Let \(N\geq\tau+\varepsilon\) be an integer. For all \(n\) such
that \(\eta_{n}(N)<y(t)\), we have \(y_{n}(\lambda_{n}(t))>0\), hence
\(\sigma_{n}\leq t\leq\tau+\varepsilon\). Therefore
\(\limsup_{n}\sigma_{n}\leq\tau+\varepsilon\) for every \(\varepsilon>0\), so that \(\limsup_{n}\sigma_{n}\leq\tau\) and \(\sigma_{n}\to\tau\).
\end{proof}

\begin{remark}[Fixed levels]\label{rem:levels}
Let \(u\in\R\), and let \(x_{n}\to x\) in \(\Dd\) and \(\Phi_{n}\to\Phi\) in
\(\cH\). Since \(\tau_{u}(x,\Phi)=\tau(x,\Phi-u)\) and
\(\Phi_{n}-u\to\Phi-u\) in \(\cH\), Theorem~\ref{thm:continuity} yields
\(\tau_{u}(x_{n},\Phi_{n})\to\tau_{u}(x,\Phi)\) whenever \((x,\Phi)\) satisfies
\NT{} at level \(u\) (Definition~\ref{def:NT}). Condition \NT{} at level \(u\) only involves the
scalarized path \(y\): if \(\tau_{u}(x,\Phi)<\infty\), it means that \(u\) is
regular for \(y\); if \(\tau_{u}(x,\Phi)=\infty\), it means that
\(\sup_{[0,T]}y<u\) for every \(T\geq 0\).
\end{remark}

In many applications, for instance in boundary value problems or in
finance, one needs the exit time jointly with the trajectory and with the
position of the path at exit, or with the path stopped at exit; the joint
law of first-passage times and positions is also a central object in the
fluctuation theory of Lévy processes (see,
e.g.,~\cite{ChaumontPellas2023,DoneyKyprianou2006,GonzalezCazaresLinMijatovic2025}).
Theorem~\ref{thm:continuity} extends to these quantities. For \(x\in\Dd\)
and \(s\in[0,\infty]\), we write \(x(\cdot\wedge s)\) for the stopped path
\(t\mapsto x(t\wedge s)\), which belongs to \(\Dd\).

\begin{corollary}[Exit position and stopped path]\label{cor:exit-position}
Let \(x_{n}\to x\) in \(\Dd\) and \(\Phi_{n}\to\Phi\) in \(\cH\), and assume
that \((x,\Phi)\) satisfies \NT{}. Let \(y(t):=\Phi(t,x(t))\),
\(\tau:=\tau(x,\Phi)\) and \(\tau^{(n)}:=\tau(x_{n},\Phi_{n})\). Then
\[
\bigl(x_{n},\Phi_{n},\tau^{(n)}\bigr)\to(x,\Phi,\tau)
\qquad\text{in } \Dd\times\cH\times[0,\infty].
\]
If, in addition, \(\tau<\infty\) and
\begin{equation*}\tag{EX}\label{eq:exit-cond}
x\ \text{is continuous at}\ \tau\qquad\text{or}\qquad y(\tau-)<0,
\end{equation*}
then \(\tau^{(n)}<\infty\) for all large \(n\), and
\[
x_{n}(\tau^{(n)})\to x(\tau)\quad\text{in } \R^{d}
\qquad\text{and}\qquad
x_{n}(\cdot\wedge\tau^{(n)})\to x(\cdot\wedge\tau)\quad\text{in } \Dd .
\]
\end{corollary}

\begin{proof}
The first convergence is a restatement of
Theorem~\ref{thm:continuity}. Assume now that \(\tau<\infty\) and that
\eqref{eq:exit-cond} holds. Since \(\tau^{(n)}\to\tau\), we have
\(\tau^{(n)}<\infty\) for all large \(n\), and we only consider such \(n\).
We use the notation of the proof of Theorem~\ref{thm:continuity}; in
particular, \(\sigma_{n}=\lambda_{n}^{-1}(\tau^{(n)})\to\tau\). Let
\[
\varepsilon_{n}:=\sup\bigl\{\|x(s)-x(\tau)\|:\ \sigma_{n}\wedge\tau\leq s\leq\sigma_{n}\vee\tau\bigr\}.
\]
We claim that \(\varepsilon_{n}\to 0\). This is clear if \(x\) is continuous
at \(\tau\). Suppose that \(y(\tau-)<0\). Then \(\tau>0\), because
\(y(0-)=y(0)\geq 0\) if \(\tau=0\); moreover, there is
\(\delta\in(0,\tau)\) such that \(y(s)\leq y(\tau-)/2\) for all
\(s\in[\tau-\delta,\tau)\), and \(\sup_{[0,\tau-\delta]}y<0\) by
\eqref{eq:NTminus}; hence \(m:=\sup_{0\leq s<\tau}y(s)<0\). Let \(N\geq\tau\)
be an integer. For all \(n\) such that \(\eta_{n}(N)<|m|\), we have
\(y_{n}(\lambda_{n}(t))\leq y(t)+\eta_{n}(N)<0\) for every \(t\in[0,\tau)\),
hence \(\sigma_{n}\geq\tau\) and
\(\varepsilon_{n}=\sup_{\tau\leq s\leq\sigma_{n}}\|x(s)-x(\tau)\|\), which
tends to \(0\) by right-continuity of \(x\) at \(\tau\).

Let \(N>\tau\) be an integer. For all \(n\) such that \(\sigma_{n}\leq N\),
\[
\begin{aligned}
\|x_{n}(\tau^{(n)})-x(\tau)\|
&\leq\|x_{n}(\lambda_{n}(\sigma_{n}))-x(\sigma_{n})\|+\|x(\sigma_{n})-x(\tau)\|\\
&\leq\sup_{0\leq s\leq N}\|x_{n}(\lambda_{n}(s))-x(s)\|+\varepsilon_{n},
\end{aligned}
\]
which tends to \(0\) by \eqref{eq:J1-char} and the claim. Moreover, since
\(\lambda_{n}\) is increasing, \(\lambda_{n}(t)\wedge\tau^{(n)}=\lambda_{n}(t\wedge\sigma_{n})\)
for all \(t\geq 0\), so that, for every \(N\in\N\),
\[
\begin{aligned}
&\sup_{0\leq t\leq N}\|x_{n}(\lambda_{n}(t)\wedge\tau^{(n)})-x(t\wedge\tau)\|\\
&\qquad\leq\sup_{0\leq s\leq N}\|x_{n}(\lambda_{n}(s))-x(s)\|
+\sup_{0\leq t\leq N}\|x(t\wedge\sigma_{n})-x(t\wedge\tau)\|\\
&\qquad\leq\sup_{0\leq s\leq N}\|x_{n}(\lambda_{n}(s))-x(s)\|+2\varepsilon_{n}\to 0 .
\end{aligned}
\]
Indeed, \(t\wedge\sigma_{n}=t\wedge\tau\) if \(t\leq\sigma_{n}\wedge\tau\), and
both \(t\wedge\sigma_{n}\) and \(t\wedge\tau\) lie in
\([\sigma_{n}\wedge\tau,\sigma_{n}\vee\tau]\) otherwise. By the
characterization \eqref{eq:J1-char}, applied with the same time changes
\(\lambda_{n}\), this gives
\(x_{n}(\cdot\wedge\tau^{(n)})\to x(\cdot\wedge\tau)\) in \(\Dd\).
\end{proof}

The first alternative in \eqref{eq:exit-cond} holds if \(x\) is continuous,
and also if \(\tau=0\), by the convention \(x(0-):=x(0)\); under it, the two
convergences are instances of the continuity of composition and stopping
maps (cf.~\cite[Chapter~13]{Whitt2002}). Since
\(y(\tau-)=\Phi(\tau,x(\tau-))\), the second alternative says that the
pre-exit position \(x(\tau-)\) lies in the open set \(\Gamma^{\Phi}_{\tau}\):
the path leaves the domain by a jump from its interior. In Route~A
(Section~\ref{sec:routeA}), the strict overshoot \(y(\tau)>0\) gives
\eqref{eq:NTplus}. For the exit-position and stopped-path conclusions at a
jump, \eqref{eq:exit-cond} requires \(y(\tau-)<0\), and this requirement
cannot simply be dropped. For instance, let
\(d=1\), \(\Phi(t,z):=z\), \(x(t):=t-1\) for \(t<1\), \(x(t):=1\) for
\(t\geq 1\), and \(x_{n}:=x+1/n\). Then \((x,\Phi)\) satisfies \NT{} with
\(\tau=1\) and \(y(\tau)=1>0\), but \(y(\tau-)=0\): the path jumps out of
the domain from its boundary. Here \(\tau^{(n)}=1-1/n\to\tau\), but
\(x_{n}(\tau^{(n)})=0\) for every \(n\), whereas \(x(\tau)=1\), and the
stopped paths
\(x_{n}(\cdot\wedge\tau^{(n)})\) converge uniformly to
\(t\mapsto\min\{t-1,0\}\), which differs from \(x(\cdot\wedge\tau)=x\).
Finally, if \(\tau=\infty\), then
\(x_{n}(\cdot\wedge\tau^{(n)})\to x=x(\cdot\wedge\tau)\) in \(\Dd\) without any
further condition: for every \(N\in\N\), we have
\(\tau^{(n)}>\lambda_{n}(N)\) for all large \(n\), so that
\(x_{n}(\lambda_{n}(t)\wedge\tau^{(n)})=x_{n}(\lambda_{n}(t))\) for all
\(t\in[0,N]\). As in Remark~\ref{rem:levels},
Corollary~\ref{cor:exit-position} applies at any level \(u\) after replacing
\(\Phi\) by \(\Phi-u\).

\subsection{Sharpness: \texorpdfstring{\NT}{(NT)} characterizes the continuity set}\label{sec:sharpness}

Theorem~\ref{thm:continuity} gives continuity at configurations satisfying
\NT{}. The next result shows that \NT{} is not merely sufficient: the
exit-time functional is discontinuous at every configuration that violates
it, even if only the barrier is perturbed, and only by constants.

\begin{proposition}[Discontinuity when \NT{} fails]\label{prop:discontinuity}
Let \((x,\Phi)\in\Dd\times\cH\) violate \NT{}. Then \(\Phi+1/n\to\Phi\) and
\(\Phi-1/n\to\Phi\) in \(\cH\), and \(\tau(x,\Phi_{n})\) does not converge to
\(\tau(x,\Phi)\) for \(\Phi_{n}:=\Phi+1/n\) if \eqref{eq:NTminus} fails, and
for \(\Phi_{n}:=\Phi-1/n\) otherwise. Consequently, the set of continuity
points of the map \((x,\Phi)\mapsto\tau(x,\Phi)\) from \(\Dd\times\cH\) to
\([0,\infty]\) is exactly the set of configurations satisfying \NT{}.
\end{proposition}

\begin{proof}
Write \(\tau:=\tau(x,\Phi)\) and \(y(t):=\Phi(t,x(t))\), and recall that
\(y(t)<0\) for \(t<\tau\).

\emph{Case 1: \eqref{eq:NTminus} fails.} Then there is a finite
\(T\in[0,\tau)\) with \(\sup_{[0,T]}y\geq 0\), and necessarily
\(\sup_{[0,T]}y=0\). Let \(\Phi_{n}:=\Phi+1/n\), so that the scalarized
path is \(y+1/n\). For each \(n\), choose \(t_{n}\in[0,T]\) with
\(y(t_{n})>-1/n\); then \(y(t_{n})+1/n>0\), hence \(\tau(x,\Phi_{n})\leq T\)
for all \(n\), while \(\tau>T\).

\emph{Case 2: \eqref{eq:NTminus} holds, \(\tau<\infty\), and
\eqref{eq:NTplus} fails.} Then there is \(\varepsilon_{0}>0\) such that
\(y(s)\leq 0\) for all \(s\in[\tau,\tau+\varepsilon_{0}]\). Let
\(\Phi_{n}:=\Phi-1/n\). The corresponding scalarized path \(y-1/n\) is at
most \(-1/n\) on \([0,\tau+\varepsilon_{0}]\), hence
\(\tau(x,\Phi_{n})\geq\tau+\varepsilon_{0}\) for all \(n\).

In both cases \(\Phi_{n}\to\Phi\) uniformly on \(E\), hence in \(\cH\), and
\(\tau(x,\Phi_{n})\not\to\tau\). Together with Theorem~\ref{thm:continuity},
this identifies the continuity set, since \(\Dd\times\cH\) is metrizable.
\end{proof}

\begin{example}[A sticking path]\label{ex:sticking}
Let \(d=1\), \(\Phi(t,z):=z\) and \(x(t):=\min\{{t-1},0\}\). Then \(y=x\),
\(\tau(x,\Phi)=1\), and \eqref{eq:NTplus} fails, since \(y\equiv 0\) on
\([1,\infty)\). For \(\Phi_{n}:=\Phi-1/n\), the scalarized path is at most
\(-1/n\) everywhere, so \(\tau(x,\Phi_{n})=\infty\) for every \(n\): an
arbitrarily small enlargement of the domain, from \((-\infty,0)\) to
\((-\infty,1/n)\), destroys the exit. For
\(\Phi_{n}:=\Phi+1/n\), we get \(\tau(x,\Phi_{n})=1-1/n\to 1\): the
convergence then holds, but only because the shift creates an overshoot.
Perturbing the path instead of the barrier, \(x_{n}:=x\pm 1/n\), has the
same effects.
\end{example}

\begin{remark}[Almost sure convergence can fail at configurations violating \NT{}]\label{rem:as-fails}
Let \(x\) and \(\Phi\) be as in Example~\ref{ex:sticking}, let
\((\beta_{n})_{n\geq 1}\) be independent Bernoulli random variables with
\(\PP(\beta_{n}=1)=1-1/n\), and set \(X_{n}:=x+(2\beta_{n}-1)/n\). Then
\(\sup_{t\geq 0}|X_{n}(t)-x(t)|=1/n\) for every \(n\) and every realization. On
\(\{\beta_{n}=1\}\) we have \(\tau(X_{n},\Phi)=1-1/n\), and on \(\{\beta_{n}=0\}\) we
have \(X_{n}\leq -1/n\), hence \(\tau(X_{n},\Phi)=\infty\). Since
\(\PP(\beta_{n}=0)=1/n\to 0\), \(\tau(X_{n},\Phi)\to 1=\tau(x,\Phi)\) in
probability. However, \(\sum_{n}\PP(\beta_{n}=0)=\infty\), so by the second
Borel--Cantelli lemma \(\tau(X_{n},\Phi)=\infty\) for infinitely many \(n\),
almost surely, and almost sure convergence fails. At configurations
violating \NT{}, the behavior of exit times thus depends on the fine
structure of the approximation.
\end{remark}

\subsection{\texorpdfstring{Weak convergence under joint convergence of \((X_{n},\Phi_{n})\)}{Weak convergence under joint convergence of (Xn, Phin)}}\label{sec:weak}

Probability spaces are assumed to be complete, and \(\Rightarrow\) denotes
convergence in distribution. A subset of a Polish space
\(\mathcal{S}\) is \emph{universally measurable} if it is measurable for the
completion of every Borel probability measure on \(\mathcal{S}\); a map from
\(\mathcal{S}\) to \([0,\infty]\) is universally measurable if the preimages of
Borel sets are universally measurable. If \(Z\) is a random element of
\(\mathcal{S}\) and \(h\) is
universally measurable, then \(h(Z)\) is a random variable whose law is
determined by the law of \(Z\).

\begin{lemma}[Measurability]\label{lem:measurability}
Let \(u\in\R\).
\begin{enumerate}[(a)]
\item The map \((x,\Phi)\mapsto\tau_{u}(x,\Phi)\) from \(\Dd\times\cH\) to
\([0,\infty]\) is universally measurable; more precisely, for every
\(a\in(0,\infty]\), the set \(\{\tau_{u}<a\}\) is an analytic subset of
\(\Dd\times\cH\).
\item The set \(\mathcal{C}_{u}\) of configurations \((x,\Phi)\) satisfying
\NT{} at level \(u\) is universally measurable.
\item The map \((x,\Phi)\mapsto\inf\{t\geq 0:\ \Phi(t,x(t))>u\}\) is Borel
measurable.
\end{enumerate}
\end{lemma}

\begin{proof}
The evaluation map \(E\times\cH\to\R\), \(((t,z),\Phi)\mapsto\Phi(t,z)\), is
continuous: if \((t_{n},z_{n})\to(t,z)\) and \(\Phi_{n}\to\Phi\) in \(\cH\),
then \(|\Phi_{n}(t_{n},z_{n})-\Phi(t,z)|\leq\sup_{K}|\Phi_{n}-\Phi|+|\Phi(t_{n},z_{n})-\Phi(t,z)|\to 0\)
for a compact set \(K\subset E\) containing all \((t_{n},z_{n})\). The map
\((x,t)\mapsto x(t)\) is Borel on \(\Dd\times[0,\infty)\): by right-continuity,
it is the pointwise limit as \(m\to\infty\) of the maps
\((x,t)\mapsto x(\lceil 2^{m}t\rceil 2^{-m})\), each of which is Borel,
because the evaluation maps \(x\mapsto x(s)\) are Borel on \(\Dd\)
(Section~\ref{sec:topologies}) and \(\{t:\ \lceil 2^{m}t\rceil=j\}\) is
Borel for every \(j\). Consequently, \((x,\Phi,t)\mapsto\Phi(t,x(t))\) is
Borel on the Polish space \(\Dd\times\cH\times[0,\infty)\), and, for every
\(t\geq 0\), \((x,\Phi)\mapsto\Phi(t,x(t))\) is Borel on \(\Dd\times\cH\).
Write \(y(t):=\Phi(t,x(t))\) and \(\tau_{u}:=\tau_{u}(x,\Phi)\).

(a) Let \(a\in(0,\infty]\) and
\(B_{a}:=\{(x,\Phi,t):\ t<a,\ \Phi(t,x(t))\geq u\}\), a Borel subset of
\(\Dd\times\cH\times[0,\infty)\). If \(\tau_{u}<a\), then \(y(\tau_{u})\geq u\)
by right-continuity, so that \((x,\Phi,\tau_{u})\in B_{a}\); conversely, if
\((x,\Phi,t)\in B_{a}\), then \(\tau_{u}\leq t<a\). Hence \(\{\tau_{u}<a\}\) is
the projection of \(B_{a}\) on \(\Dd\times\cH\). Projections of Borel subsets
of Polish spaces are analytic, and analytic sets are universally measurable~\cite[Sections~14.A and 29.B]{Kechris1995}. Since the sets \([0,a)\)
generate the Borel \(\sigma\)-field of \([0,\infty]\), (a) follows.

(b) For a rational \(q\geq 0\), right-continuity gives
\(\sup_{[0,q]}y=\sup_{r\in\Q\cap[0,q]}y(r)\), which is a Borel function of
\((x,\Phi)\). If \(T<\tau_{u}\) is finite, there is a rational
\(q\in[T,\tau_{u})\), and \(\sup_{[0,T]}y\leq\sup_{[0,q]}y\). Hence
\((x,\Phi)\) satisfies \eqref{eq:NTminus} at level \(u\) if and only if, for
every rational \(q\geq 0\), either \(\tau_{u}\leq q\) or
\(\sup_{r\in\Q\cap[0,q]}y(r)<u\). Similarly, by right-continuity,
\eqref{eq:NTplus} holds at level \(u\) if and only if \(\tau_{u}=\infty\) or,
for every \(k\in\N\), there is a rational \(r\) with
\(\tau_{u}\leq r<\tau_{u}+1/k\) and \(y(r)>u\). These are countable Boolean
combinations of universally measurable sets, by (a), and of Borel sets.

(c) Let \(a\in(0,\infty]\). If \(y(t)>u\) for some \(t<a\), then \(y(r)>u\) for all rationals
\(r\in[t,a)\) close enough to \(t\), by right-continuity. Hence
\[
\bigl\{\inf\{t\geq 0:\ y(t)>u\}<a\bigr\}=\bigcup_{r\in\Q\cap[0,a)}\{y(r)>u\},
\]
which is Borel.
\end{proof}

\begin{remark}[On Borel measurability]\label{rem:measurability}
We do not claim that \(\tau_{u}\) is Borel measurable. The argument used for
Lemma~\ref{lem:measurability}(c), which concerns first passage into the
open set \((u,\infty)\), does not apply to first passage into the closed set
\([u,\infty)\): since the supremum of a càdlàg function over a compact
interval need not be attained as a value (Remark~\ref{rem:NTminus}), the
inequality \(\sup_{r\in\Q\cap[0,q]}y(r)\geq u\) does not imply that
\(y(t)\geq u\) for some \(t\leq q\). By contrast, the set \(\mathcal{C}_{u}\)
of Lemma~\ref{lem:measurability}(b) is even Borel. Indeed, by
Proposition~\ref{prop:discontinuity}, applied to \(\Phi-u\), and since
\((x,\Phi)\mapsto(x,\Phi-u)\) is a homeomorphism of \(\Dd\times\cH\),
\(\mathcal{C}_{u}\) is the set of continuity points of
\(\tau_{u}:\Dd\times\cH\to[0,\infty]\); and the set of continuity points of
any map between metrizable spaces is a countable intersection of open sets,
namely of the sets on which its oscillation is smaller than \(1/k\),
\(k\in\N\). This does not give the Borel measurability of \(\tau_{u}\)
itself. Universal measurability is all that is needed below.
\end{remark}

\begin{theorem}[Weak convergence of exit times]\label{thm:weak}
Let \((X_{n},\Phi_{n})\Rightarrow(X,\Phi)\) in \(\Dd\times\cH\), and assume
that
\begin{equation}\label{eq:NT-as}
\PP\bigl((X,\Phi)\text{ satisfies \NT{}}\bigr)=1 .
\end{equation}
Then \(\tau(X_{n},\Phi_{n})\Rightarrow\tau(X,\Phi)\) in \([0,\infty]\).
\end{theorem}

\begin{proof}
By Lemma~\ref{lem:measurability}, \(\tau^{(n)}:=\tau(X_{n},\Phi_{n})\) and
\(\tau:=\tau(X,\Phi)\) are random variables, their laws are determined by
those of \((X_{n},\Phi_{n})\) and \((X,\Phi)\), and the probability in
\eqref{eq:NT-as} is well defined. Since \(\Dd\times\cH\) is Polish,
Skorokhod's representation theorem~\cite[Theorem~6.7]{Billingsley1999}
provides random elements \((\tilde{X}_{n},\tilde{\Phi}_{n})\) and
\((\tilde{X},\tilde{\Phi})\), defined on a common complete probability
space, with the same laws as \((X_{n},\Phi_{n})\) and \((X,\Phi)\), such that
\((\tilde{X}_{n},\tilde{\Phi}_{n})\to(\tilde{X},\tilde{\Phi})\) almost surely.
By \eqref{eq:NT-as}, \((\tilde{X},\tilde{\Phi})\) satisfies \NT{} almost
surely, and Theorem~\ref{thm:continuity} yields
\(\tau(\tilde{X}_{n},\tilde{\Phi}_{n})\to\tau(\tilde{X},\tilde{\Phi})\) almost
surely. Almost sure convergence implies convergence in distribution, and
\(\tau(\tilde{X}_{n},\tilde{\Phi}_{n})\) and \(\tau(\tilde{X},\tilde{\Phi})\) have
the same laws as \(\tau^{(n)}\) and \(\tau\).
\end{proof}

\begin{remark}[Dependence between \(X\) and \(\Phi\)]\label{rem:dependence}
No independence or adaptedness is assumed: \(\Phi\) may be random and
coupled to \(X\). The argument is pathwise once the joint convergence
\((X_{n},\Phi_{n})\Rightarrow(X,\Phi)\) and condition \NT{} for the limit
\((X,\Phi)\), which only involves the scalarized limit \(t\mapsto\Phi(t,X(t))\)
(Remark~\ref{rem:levels}), are established. In the example of
Section~\ref{sec:example-coupled}, where the approximating pairs are denoted
by \((X^{(n)},\Phi^{(n)})\), the barrier \(\Phi^{(n)}\) is built
from the path \(X^{(n)}\) and its increments, need not be adapted to the
natural filtration of \(X^{(n)}\) (it is not adapted as soon as the squared
norms of the increments are not almost surely constant), and converges to a
random barrier that is a functional of the limit path.
\end{remark}

Corollary~\ref{cor:exit-position} has the following counterpart in
distribution.

\begin{corollary}[Joint convergence, exit position and stopped path]\label{cor:joint}
Let \((X_{n},\Phi_{n})\Rightarrow(X,\Phi)\) in \(\Dd\times\cH\), assume
\eqref{eq:NT-as}, and let \(\tau^{(n)}:=\tau(X_{n},\Phi_{n})\) and
\(\tau:=\tau(X,\Phi)\). Then
\begin{equation}\label{eq:joint-weak}
\bigl(X_{n},\Phi_{n},\tau^{(n)}\bigr)\Rightarrow(X,\Phi,\tau)
\qquad\text{in } \Dd\times\cH\times[0,\infty].
\end{equation}
If, in addition, almost surely, \(\tau<\infty\) and \((X,\Phi)\) satisfies
\eqref{eq:exit-cond}, that is, \(X\) is continuous at \(\tau\) or
\(\Phi(\tau,X(\tau-))<0\), then \(\PP(\tau^{(n)}<\infty)\to 1\) and
\begin{equation}\label{eq:joint-exit}
\bigl(X_{n},\Phi_{n},\tau^{(n)},X_{n}(\tau^{(n)}),X_{n}(\cdot\wedge\tau^{(n)})\bigr)
\Rightarrow\bigl(X,\Phi,\tau,X(\tau),X(\cdot\wedge\tau)\bigr)
\end{equation}
in \(\Dd\times\cH\times[0,\infty]\times\R^{d}\times\Dd\), where
\(X_{n}(\tau^{(n)}):=0\) on \(\{\tau^{(n)}=\infty\}\).
\end{corollary}

\begin{proof}
We first check measurability. The map \((x,t)\mapsto x(t)\) is Borel on
\(\Dd\times[0,\infty)\) (proof of Lemma~\ref{lem:measurability}), and so is
\((x,t)\mapsto x(t-)\), as the pointwise limit, as \(m\to\infty\), of the
Borel maps \((x,t)\mapsto x\bigl((\lceil 2^{m}t\rceil-1)^{+}2^{-m}\bigr)\).
Since the Borel \(\sigma\)-field of \(\Dd\) is generated by the evaluation
maps (Section~\ref{sec:topologies}), the map
\((x,s)\mapsto x(\cdot\wedge s)\) is Borel from \(\Dd\times[0,\infty]\) to
\(\Dd\). Together with Lemma~\ref{lem:measurability}(a) and the continuity
of the evaluation map \(((t,z),\Phi)\mapsto\Phi(t,z)\), this shows that the
map
\[
\Xi(x,\Phi):=\bigl(x,\Phi,\tau(x,\Phi),x(\tau(x,\Phi)),x(\cdot\wedge\tau(x,\Phi))\bigr),
\]
where \(x(\tau(x,\Phi)):=0\) if \(\tau(x,\Phi)=\infty\), is universally
measurable, and so is the set \(\mathcal{C}'\) of configurations
\((x,\Phi)\) such that \(\tau(x,\Phi)<\infty\) and \eqref{eq:exit-cond}
holds.

Let \((\tilde{X}_{n},\tilde{\Phi}_{n})\to(\tilde{X},\tilde{\Phi})\) almost
surely be as in the proof of Theorem~\ref{thm:weak}. By
Theorem~\ref{thm:continuity}, the first three coordinates of
\(\Xi(\tilde{X}_{n},\tilde{\Phi}_{n})\) converge almost surely to those of
\(\Xi(\tilde{X},\tilde{\Phi})\). Under the additional assumption,
\((\tilde{X},\tilde{\Phi})\in\mathcal{C}'\) almost surely, and
Corollary~\ref{cor:exit-position} gives, almost surely,
\(\tau(\tilde{X}_{n},\tilde{\Phi}_{n})<\infty\) for all large \(n\) and
\(\Xi(\tilde{X}_{n},\tilde{\Phi}_{n})\to\Xi(\tilde{X},\tilde{\Phi})\). Since
\(\Xi(\tilde{X}_{n},\tilde{\Phi}_{n})\) and \(\Xi(\tilde{X},\tilde{\Phi})\) have
the same laws as \(\Xi(X_{n},\Phi_{n})\) and \(\Xi(X,\Phi)\), this proves
\eqref{eq:joint-weak}, \(\PP(\tau^{(n)}<\infty)\to 1\) and
\eqref{eq:joint-exit}.
\end{proof}

The additional assumption holds whenever \(X\) has continuous paths and
\(\tau<\infty\) almost surely, as in Section~\ref{sec:example-coupled} and,
when the boundary \(g\) is bounded above, in
Section~\ref{sec:example-donsker}. If \(\tau=\infty\) with positive
probability, the stopped paths still converge,
\(X_{n}(\cdot\wedge\tau^{(n)})\Rightarrow X(\cdot\wedge\tau)\), jointly with
\eqref{eq:joint-weak}, as soon as \eqref{eq:exit-cond} holds almost surely
on \(\{\tau<\infty\}\), by the observation on the case \(\tau=\infty\) that
follows Corollary~\ref{cor:exit-position}.

\section{Exit-time profiles and functional convergence in the \texorpdfstring{\(M_{1}\)}{M1} topology}\label{sec:profiles}

This section establishes a functional limit theorem for exit-time
profiles indexed by the barrier level \(u\). The profile
\(u\mapsto\tau_{u}\) is monotone and typically has jumps that are not well
aligned under approximation. Hence \(J_{1}\) is in general too strong,
whereas the \(M_{1}\) topology is natural for monotone càdlàg functions~\cite{PuhalskiiWhitt1997,Whitt2002}.

\subsection{Nested domains and exit-time profiles}\label{sec:profile-setup}

Let \(\Phi\in\cH\), \(x\in\Dd\) and \(y(t):=\Phi(t,x(t))\). The domains
\(\Gamma^{\Phi,u}\), \(u\in\R\), defined in \eqref{eq:barrier-domain}, are
nested, and the exit times \(\tau_{u}(x,\Phi)=T_{y}(u)\) are nondecreasing
in \(u\). By Lemma~\ref{lem:regular-continuity}, the map
\(u\mapsto T_{y}(u)\) is in general neither left- nor right-continuous, and
its right-continuous version is
\[
T^{+}_{y}(u)=\inf\{t\geq 0:\ \Phi(t,x(t))>u\},
\]
the first time at which the path leaves the closed sets
\(\{z:\ \Phi(t,z)\leq u\}\). The two functions coincide at every level at which \((x,\Phi)\) satisfies \NT{}: at
regular levels by Lemma~\ref{lem:regular-continuity}, and trivially at
levels \(u\) with \(T_{y}(u)=\infty\). We work with \(T^{+}_{y}\).

Fix \(u_{0}<u_{1}\). Profiles may take the value \(+\infty\) at levels at
which no exit occurs. We therefore let \(\Ms[u_{0},u_{1}]\) be the set of
nondecreasing, right-continuous functions \(f:[u_{0},u_{1}]\to[0,\infty]\),
and we let \(\gbar:[0,\infty]\to[0,1]\), \(\gbar(t):=t/(1+t)\),
\(\gbar(\infty):=1\). For \(a<b\), let \(D[a,b]\) denote the space of
real functions on \([a,b]\) that are right-continuous on \([a,b)\) and have
left limits on \((a,b]\); unless another topology is specified, it is
equipped with the Skorokhod \(M_{1}\) topology (see
\hyperref[app:M1]{Appendix~A} and~\cite[Chapter~12]{Whitt2002}). For
\(f\in\Ms[u_{0},u_{1}]\), the function \(\gbar\circ f\) is a nondecreasing
element of \(D[u_{0},u_{1}]\) with values in \([0,1]\), and we equip
\(\Ms[u_{0},u_{1}]\) with the metric
\[
d_{\Ms}(f,h):=d_{M_{1}}(\gbar\circ f,\gbar\circ h),
\]
where \(d_{M_{1}}\) is the \(M_{1}\) metric \eqref{eq:M1-metric}. We refer to
convergence in \(\Ms[u_{0},u_{1}]\) as \(M_{1}\) convergence of profiles.
This is consistent with the usual \(M_{1}\) topology: if \(f_{n}\) and \(f\)
are finite-valued elements of \(\Ms[u_{0},u_{1}]\), then \(f_{n}\to f\) in
\(\Ms[u_{0},u_{1}]\) if and only if \(f_{n}\to f\) in \((D[u_{0},u_{1}],M_{1})\).
Indeed, \(\gbar\) is an increasing bijection from \([0,\infty)\) onto \([0,1)\)
which is \(1\)-Lipschitz and whose inverse is Lipschitz on every interval
\([0,c]\), \(c<1\); composition with \(\gbar\) maps the parametric
representations of \(f\) (\hyperref[app:M1]{Appendix~A}) onto those of
\(\gbar\circ f\); and \(M_{1}\) convergence implies convergence at the right
end point, so that, under either mode of convergence, \(f_{n}(u_{1})\to
f(u_{1})<\infty\) and the functions \(f_{n}\) are eventually bounded by a
common constant.

The unregularized profile \(T_{y}\), that is, \(u\mapsto\tau_{u}(x,\Phi)\),
may differ from \(T^{+}_{y}\) at individual levels, but the two profiles
have the same completed graph on \([u_{0},u_{1}]\) whenever they agree at
\(u_{0}\) and \(u_{1}\), in particular if \((x,\Phi)\) satisfies \NT{} at
these two levels. Here, for a nondecreasing function \(g\) on
\([u_{0},u_{1}]\) that need not be right-continuous, the completed graph is
formed by the vertical segments joining \(g(u-)\) to \(g(u+)\) at each
interior level \(u\), \(g(u_{0})\) to \(g(u_{0}+)\) at \(u_{0}\), and
\(g(u_{1}-)\) to \(g(u_{1})\) at \(u_{1}\); for right-continuous \(g\), this is
the completed graph of \hyperref[app:M1]{Appendix~A}, on which the metric
\(d_{\Ms}\) is based after composition with \(\gbar\). The claim holds
because \(T_{y}(v)\leq T^{+}_{y}(v)\leq T_{y}(v')\) for all \(v<v'\), so that
\(T_{y}\) and \(T^{+}_{y}\) have the same left limits and the same right
limits at every level. The right-continuous profile \(T^{+}_{y}\) remains
the one that belongs to \(\Ms[u_{0},u_{1}]\).

\subsection{Main theorem: functional convergence of exit-time profiles}\label{sec:profile-thm}

Let \((X_{n},\Phi_{n})\) and \((X,\Phi)\) be random elements of
\(\Dd\times\cH\), let \(Y_{n}(t):=\Phi_{n}(t,X_{n}(t))\) and
\(Y(t):=\Phi(t,X(t))\) be the scalarized processes, and define the
exit-time profiles
\begin{equation}\label{eq:profiles}
\tau^{+}_{n}(u):=T^{+}_{Y_{n}}(u),\qquad\tau^{+}(u):=T^{+}_{Y}(u),\qquad u\in[u_{0},u_{1}].
\end{equation}
Thus \(\tau^{+}(u)=\tau_{u}(X,\Phi)\) at every level \(u\) at which
\((X,\Phi)\) satisfies \NT{}, and similarly for \(\tau^{+}_{n}\). By
Lemma~\ref{lem:measurability}(c), \(\tau^{+}_{n}(u)\) and \(\tau^{+}(u)\) are
random variables for every \(u\). Since the Borel \(\sigma\)-field of
\((D[u_{0},u_{1}],M_{1})\) is generated by the evaluation maps~\cite[Chapters~11 and~12]{Whitt2002}, \(\tau^{+}_{n}\) and \(\tau^{+}\) are random
elements of \(\Ms[u_{0},u_{1}]\).

\begin{theorem}[Functional convergence of exit-time profiles in \(M_{1}\)]\label{thm:profiles}
Let \({u_{0}<u_{1}}\). Assume that \((X_{n},\Phi_{n})\Rightarrow(X,\Phi)\) in
\(\Dd\times\cH\) and that, almost surely, \((X,\Phi)\) satisfies \NT{} at the
levels \(u_{0}\) and \(u_{1}\) (Definition~\ref{def:NT}). Then
\[
\tau^{+}_{n}\Rightarrow\tau^{+}\qquad\text{in } \Ms[u_{0},u_{1}].
\]
\end{theorem}

The proof reduces the problem to the almost sure \(J_{1}\) convergence of
the scalarized paths, proves pointwise convergence of the profiles at the
levels at which \NT{} holds, and upgrades pointwise convergence on a dense
set of levels to \(M_{1}\) convergence by monotonicity. It relies on the
following three lemmas. Here and below, a path \(y\in\Done\) is said to
satisfy \NT{} at level \(u\) if the configuration \((y,\iota)\) does, where
\(d=1\) and \(\iota(t,z):=z\) is the identity barrier.

\begin{lemma}[Pointwise convergence at a level satisfying \NT{}]\label{lem:pointwise}
Let \(y_{n}\to y\) in \(\Done\), and let \(u\in\R\) be a level at which \(y\)
satisfies \NT{}. Then \(T_{y_{n}}(u)\to T_{y}(u)\) and
\(T^{+}_{y_{n}}(u)\to T^{+}_{y}(u)=T_{y}(u)\) in \([0,\infty]\).
\end{lemma}

\begin{proof}
The convergence \(T_{y_{n}}(u)\to T_{y}(u)\) is Theorem~\ref{thm:continuity}
(with Remark~\ref{rem:levels}) applied with \(d=1\), \(x_{n}:=y_{n}\),
\(x:=y\) and \(\Phi_{n}=\Phi:=\iota\). Moreover,
\(T^{+}_{y}(u)=T_{y}(u)\), as observed in Section~\ref{sec:profile-setup}. If
\(T_{y}(u)=\infty\), then \(T^{+}_{y_{n}}(u)\geq T_{y_{n}}(u)\to\infty\). If
\(T_{y}(u)<\infty\), it remains to show that
\(\limsup_{n}T^{+}_{y_{n}}(u)\leq T_{y}(u)\), because
\(T^{+}_{y_{n}}(u)\geq T_{y_{n}}(u)\). Let \(\lambda_{n}\in\Lambda\) satisfy
\eqref{eq:J1-char} for \(y_{n}\to y\), and let \(\varepsilon>0\). By
\eqref{eq:Rplus}, there is \(t\in[T_{y}(u),T_{y}(u)+\varepsilon]\) with
\(y(t)>u\). For \(n\) large enough, \(y_{n}(\lambda_{n}(t))>u\), and hence
\(T^{+}_{y_{n}}(u)\leq\lambda_{n}(t)\leq T_{y}(u)+\varepsilon+\sup_{s\geq 0}|\lambda_{n}(s)-s|\).
\end{proof}

\begin{lemma}[Monotone \(M_{1}\) criterion]\label{lem:monotone-M1}
Let \(a<b\), and let \(f_{n},f\in D[a,b]\) be nondecreasing. If
\(f_{n}(t)\to f(t)\) for every \(t\) in a dense subset of \([a,b]\) that
contains \(a\) and \(b\), then \(f_{n}\to f\) in the \(M_{1}\) topology.
\end{lemma}

This criterion is classical (see~\cite[Chapter~12]{Whitt2002}); a self-contained proof
is given in \hyperref[app:M1]{Appendix~A}.

\begin{lemma}[Levels violating \NT{} are countable]\label{lem:countable}
For every \(y\in\Done\), the set \(\mathcal{B}(y)\) of levels \(u\in\R\) at
which \(y\) does not satisfy \NT{} is countable.
\end{lemma}

\begin{proof}
If \(u\in\mathcal{B}(y)\) and \(T_{y}(u)<\infty\), then \(u\) is not regular
for \(y\), so that the nondecreasing function \(T_{y}\) is discontinuous at
\(u\) by Lemma~\ref{lem:regular-continuity}; there are at most countably
many such levels. If \(u\in\mathcal{B}(y)\) and \(T_{y}(u)=\infty\), then
\(y(t)<u\) for all \(t\geq 0\), whereas \(\sup_{[0,T]}y\geq u\) for some
\(T\geq 0\); hence \(u=\sup_{t\geq 0}y(t)\), and there is at most one such
level.
\end{proof}

\begin{proof}[Proof of Theorem~\ref{thm:profiles}]
By Lemma~\ref{lem:scalarization}, the scalarization map \(S\) is
continuous, so the continuous mapping theorem~\cite[Theorem~2.7]{Billingsley1999} yields \(Y_{n}\Rightarrow Y\) in
\(\Done\). Since \(\Done\) is Polish, Skorokhod's representation theorem~\cite[Theorem~6.7]{Billingsley1999} provides random elements
\(\tilde{Y}_{n}\) and \(\tilde{Y}\) of \(\Done\), defined on a common complete
probability space, with the same laws as \(Y_{n}\) and \(Y\), such that
\(\tilde{Y}_{n}\to\tilde{Y}\) almost surely. Condition \NT{} at level \(u\) only
depends on the scalarized path (Remark~\ref{rem:levels}), and the set of
paths satisfying it is universally measurable
(Lemma~\ref{lem:measurability}(b)); hence, almost surely, \(\tilde{Y}\)
satisfies \NT{} at the levels~\(u_{0}\) and~\(u_{1}\).

Fix an outcome for which \(\tilde{Y}_{n}\to\tilde{Y}\) and
\(u_{0},u_{1}\notin\mathcal{B}(\tilde{Y})\), and let
\(D_{0}:=[u_{0},u_{1}]\setminus\mathcal{B}(\tilde{Y})\). By
Lemma~\ref{lem:countable}, \(D_{0}\) is dense in \([u_{0},u_{1}]\), and it
contains \(u_{0}\) and \(u_{1}\). By Lemma~\ref{lem:pointwise},
\(\gbar(T^{+}_{\tilde{Y}_{n}}(u))\to\gbar(T^{+}_{\tilde{Y}}(u))\) for every
\(u\in D_{0}\). The functions \(\gbar\circ T^{+}_{\tilde{Y}_{n}}\) and
\(\gbar\circ T^{+}_{\tilde{Y}}\) are nondecreasing elements of
\(D[u_{0},u_{1}]\), so Lemma~\ref{lem:monotone-M1} yields
\(\gbar\circ T^{+}_{\tilde{Y}_{n}}\to\gbar\circ T^{+}_{\tilde{Y}}\) in \(M_{1}\), that
is, \(T^{+}_{\tilde{Y}_{n}}\to T^{+}_{\tilde{Y}}\) in \(\Ms[u_{0},u_{1}]\).

Thus the profiles of \(\tilde{Y}_{n}\) converge almost surely, hence in
distribution, to the profile of \(\tilde{Y}\). Since the profile map
\(y\mapsto T^{+}_{y}|_{[u_{0},u_{1}]}\) is Borel measurable from \(\Done\) to
\(\Ms[u_{0},u_{1}]\) (see the discussion preceding the theorem), these
profiles have the same laws as \(\tau^{+}_{n}\) and \(\tau^{+}\), which
proves the theorem.
\end{proof}

The same arguments give the convergence of the finite-dimensional
distributions of the profiles at fixed levels satisfying \NT{}.

\begin{corollary}[Finite-dimensional convergence]\label{cor:fdd}
Let \((X_{n},\Phi_{n})\Rightarrow(X,\Phi)\) in \(\Dd\times\cH\), and let
\(v_{1}<\cdots<v_{m}\) be deterministic levels at each of which \((X,\Phi)\)
satisfies \NT{} almost surely. Then
\begin{equation}\label{eq:fdd}
\bigl(\tau_{v_{i}}(X_{n},\Phi_{n})\bigr)_{1\leq i\leq m}
\Rightarrow\bigl(\tau_{v_{i}}(X,\Phi)\bigr)_{1\leq i\leq m}
\qquad\text{in } [0,\infty]^{m}.
\end{equation}
The same convergence holds, with the same limit, if
\(\tau_{v_{i}}(X_{n},\Phi_{n})\) is replaced by \(T^{+}_{Y_{n}}(v_{i})\) for
\(1\leq i\leq m\), and both convergences hold jointly with each other and
with \((X_{n},\Phi_{n})\Rightarrow(X,\Phi)\).
\end{corollary}

\begin{proof}
Let \((\tilde{X}_{n},\tilde{\Phi}_{n})\to(\tilde{X},\tilde{\Phi})\) almost
surely be as in the proof of Theorem~\ref{thm:weak}, with scalarized paths
\(\tilde{Y}_{n}\) and \(\tilde{Y}\), so that \(\tilde{Y}_{n}\to\tilde{Y}\) in
\(\Done\) almost surely by Lemma~\ref{lem:scalarization}. Almost surely,
\((\tilde{X},\tilde{\Phi})\), and hence \(\tilde{Y}\), satisfies \NT{} at each
of the levels \(v_{1},\ldots,v_{m}\) (Lemma~\ref{lem:measurability}(b) and
Remark~\ref{rem:levels}), and Lemma~\ref{lem:pointwise} then gives
\(T_{\tilde{Y}_{n}}(v_{i})\to T_{\tilde{Y}}(v_{i})\) and
\(T^{+}_{\tilde{Y}_{n}}(v_{i})\to T^{+}_{\tilde{Y}}(v_{i})=T_{\tilde{Y}}(v_{i})\)
for every \(i\). By Lemma~\ref{lem:measurability}(a) and~(c),
\(T_{\tilde{Y}_{n}}(v_{i})=\tau_{v_{i}}(\tilde{X}_{n},\tilde{\Phi}_{n})\) is a
universally measurable function of \((\tilde{X}_{n},\tilde{\Phi}_{n})\) and
\(T^{+}_{\tilde{Y}_{n}}(v_{i})\) a Borel one, and similarly in the limit.
Hence the joint laws are those of the corresponding quantities without
tildes, and the claims follow.
\end{proof}

\subsection{Remarks and examples}\label{sec:profile-remarks}

\begin{remark}[Finite profiles]\label{rem:finite-profiles}
If, in addition, \(\tau^{+}(u_{1})<\infty\) almost surely, then \(\tau^{+}\)
is almost surely a finite-valued element of \(D[u_{0},u_{1}]\),
\({\PP(\tau^{+}_{n}(u_{1})<\infty)}\to 1\), and
\(\E[h(\tau^{+}_{n})\mathbf{1}_{\{\tau^{+}_{n}(u_{1})<\infty\}}]\to\E[h(\tau^{+})]\)
for every bounded function \(h:D[u_{0},u_{1}]\to\R\) that is continuous for
\(M_{1}\), with the convention that the product inside the expectation is
\(0\) on \(\{\tau^{+}_{n}(u_{1})=\infty\}\).
Indeed, the set \(O:=\{f\in\Ms[u_{0},u_{1}]:\ f(u_{1})<\infty\}\) is open in
\(\Ms[u_{0},u_{1}]\), the topology of \(\Ms[u_{0},u_{1}]\) coincides with the
\(M_{1}\) topology on \(O\) (Section~\ref{sec:profile-setup}), and it
suffices to apply the portmanteau and continuous mapping theorems to the
function \(f\mapsto h(f)\mathbf{1}_{O}(f)\), which is continuous at every
point of \(O\).
\end{remark}

\begin{remark}[Why \NT{} is only needed at the end levels]\label{rem:end-levels}
Since profiles are monotone, \(M_{1}\) convergence follows from pointwise
convergence on any dense set of levels containing the end levels
(Lemma~\ref{lem:monotone-M1}). For a fixed limit path, \NT{} fails at most at
countably many levels, namely at levels \(u\) with \(T_{Y}(u)<\infty\) at
which the function \(T_{Y}\) is discontinuous, and possibly at the level
\(\sup_{t}Y(t)\) (Lemma~\ref{lem:countable}); pointwise convergence holds at
all other levels (Lemma~\ref{lem:pointwise}). These exceptional levels are random,
but no condition on them is needed, because the dense set in
Lemma~\ref{lem:monotone-M1} may depend on the realization. Only the end
levels require a separate condition. In the examples of
Section~\ref{sec:verification}, \NT{} is verified at each fixed
deterministic level, almost surely: at every fixed \(u\in\R\) in
Section~\ref{sec:example-donsker}, at every fixed \(u>-\rho_{*}^{2}\) in
Section~\ref{sec:example-coupled}, where \(\rho_{*}>0\) is a lower bound for
the tolerance \(\rho\), and at every fixed \(u>-h\) in
Remark~\ref{rem:drawdown}, where \(h>0\) is the drawdown threshold. Applied to \(u_{0}\) and \(u_{1}\), this is all that
Theorem~\ref{thm:profiles} requires. It does not mean that, almost surely,
\NT{} holds simultaneously at all these levels, and in general it does not:
if the scalarized path \(Y\) is continuous and attains its maximum over
\([0,1]\) at some time strictly before \(1\), then \eqref{eq:Rplus} fails at
the random level \(m:=\max_{0\leq t\leq 1}Y(t)\), because \(Y(t)\leq m\) for
\(T_{Y}(m)\leq t\leq 1\). This happens with positive probability in
Section~\ref{sec:example-donsker} and in Section~\ref{sec:example-coupled}
with constant tolerance \(\rho\).
\end{remark}

\begin{remark}[Condition \NT{} at the end levels cannot be dropped]\label{rem:end-sharp}
Fix \({u^{*}\in\R}\) and consider the path
\[
y(t):=
\begin{cases}
u^{*}-(1-t)^{2}, & 0\leq t<1,\\
u^{*}, & 1\leq t<2,\\
u^{*}+1, & t\geq 2,
\end{cases}
\]
which reaches \(u^{*}\) at time \(1\) but does not exceed it before time
\(2\), so that \eqref{eq:Rplus} fails at \(u^{*}\). Let
\(y_{n}:=y+1/n\), so that \(y_{n}\to y\) uniformly. Then
\(T^{+}_{y_{n}}(u)=T^{+}_{y}(u-1/n)\) for all \(u\), and in particular
\[
T^{+}_{y_{n}}(u^{*})=1-n^{-1/2}\to 1,\qquad\text{whereas}\qquad T^{+}_{y}(u^{*})=2 .
\]
On \([u_{0},u_{1}]:=[u^{*}-1/4,u^{*}]\), \NT{} holds at \(u_{0}\) but fails at
\(u_{1}=u^{*}\), and the profiles of \(y_{n}\) do not converge to that of
\(y\) in \(M_{1}\), because \(M_{1}\) convergence implies convergence at the
right end point. By contrast, on any interval \([u_{0},u_{1}]\) that
contains \(u^{*}\) in its interior and at whose end levels \NT{} holds, the
profiles do converge in \(M_{1}\), although they do not converge at the level~\(u^{*}\).
\end{remark}

\begin{remark}[Open level intervals]\label{rem:open-intervals}
Theorem~\ref{thm:profiles} extends to open, possibly unbounded, level
intervals. Let \(I\subset\R\) be an open interval, and let
\((u_{0,k})_{k\geq 1}\) and \((u_{1,k})_{k\geq 1}\) be sequences in \(I\) such
that \(u_{0,k}<u_{1,k}\), \(u_{0,k}\downarrow\inf I\) and
\(u_{1,k}\uparrow\sup I\), and such that \((X,\Phi)\) satisfies \NT{} almost
surely at each of the levels \(u_{0,k}\) and \(u_{1,k}\). On the probability
space used in the proof of Theorem~\ref{thm:profiles}, almost surely,
\(\tilde{Y}\) then satisfies \NT{} at all these countably many levels, and
the pathwise part of that proof gives
\(T^{+}_{\tilde{Y}_{n}}\to T^{+}_{\tilde{Y}}\) in \(\Ms[u_{0,k},u_{1,k}]\) for
every \(k\geq 1\) simultaneously. As in that proof, each restriction, and
hence the sequence of restrictions, is a Borel function of the scalarized
path, so that the laws transfer and the restrictions converge jointly in
distribution:
\[
\bigl(T^{+}_{Y_{n}}|_{[u_{0,k},u_{1,k}]}\bigr)_{k\geq 1}\Rightarrow
\bigl(T^{+}_{Y}|_{[u_{0,k},u_{1,k}]}\bigr)_{k\geq 1}
\qquad\text{in } \prod_{k\geq 1}\Ms[u_{0,k},u_{1,k}],
\]
where the product carries the product topology. By \NT{}, the limit profile
\(T^{+}_{Y}\) is almost surely continuous, as a \([0,\infty]\)-valued
function, at every end level \(u\in\{u_{0,k},u_{1,k}\}\): by
Lemma~\ref{lem:regular-continuity} if \(T_{Y}(u)<\infty\), and because
\eqref{eq:NTminus} at level \(u\) forces \(T_{Y}(v)\to\infty\) as
\(v\uparrow u\) if \(T_{Y}(u)=\infty\). The above convergence therefore
agrees with the usual notion of \(M_{1}\) convergence on a noncompact
interval, which is defined through restrictions to compact subintervals
whose end points are continuity points of the limit
(cf.~\cite[Section~12.9]{Whitt2002}). One may take \(I=\R\) in
Section~\ref{sec:example-donsker}, \(I=(-\rho_{*}^{2},\infty)\) in
Section~\ref{sec:example-coupled}, and \(I=(-h,\infty)\) in
Remark~\ref{rem:drawdown}.
\end{remark}

The next example shows that \(M_{1}\) cannot be replaced by \(J_{1}\) in
Theorem~\ref{thm:profiles}: a single jump of the limit profile may be
approximated by a cluster of small jumps, which \(M_{1}\) regards as
convergent whereas \(J_{1}\) does not.

\begin{example}[\(J_{1}\) convergence of profiles can fail]\label{ex:J1-fails}
Let \([u_{0},u_{1}]:=[1/4,3/4]\), let
\[
y(t):=
\begin{cases}
0, & t<1,\\
1/2, & 1\leq t<2,\\
1, & t\geq 2,
\end{cases}
\]
and, for \(n\geq 3\), let
\[
y_{n}(t):=
\begin{cases}
0, & t<1,\\
\dfrac{1}{2}+\dfrac{k}{n^{2}}, & t\in\Bigl[1+\dfrac{k}{n},1+\dfrac{k+1}{n}\Bigr),\quad k=0,\ldots,n-1,\\[1ex]
1, & t\geq 2.
\end{cases}
\]
Then \(\sup_{t\geq 0}|y_{n}(t)-y(t)|<1/n\), so \(y_{n}\to y\) uniformly, and in
particular in \(\Done\). The profiles on \([u_{0},u_{1}]\) are
\[
T^{+}_{y}(u)=
\begin{cases}
1, & u\in[1/4,1/2),\\
2, & u\in[1/2,3/4],
\end{cases}
\]
and
\[
T^{+}_{y_{n}}(u)=
\begin{cases}
1, & u\in[1/4,1/2),\\
1+\dfrac{k}{n}, & u\in\Bigl[\dfrac{1}{2}+\dfrac{k-1}{n^{2}},\dfrac{1}{2}+\dfrac{k}{n^{2}}\Bigr),\quad k=1,\ldots,n-1,\\[1ex]
2, & u\in\Bigl[\dfrac{1}{2}+\dfrac{n-1}{n^{2}},\dfrac{3}{4}\Bigr].
\end{cases}
\]
The path \(y\) satisfies \NT{} at every level of \([1/4,3/4]\) except
\(u=1/2\), at which it sticks. Accordingly, \(T^{+}_{y_{n}}(u)\to T^{+}_{y}(u)\)
for every \(u\neq 1/2\), and \(T^{+}_{y_{n}}\to T^{+}_{y}\) in \(M_{1}\) by
Lemma~\ref{lem:monotone-M1}. However, \(T^{+}_{y_{n}}\) does not converge to
\(T^{+}_{y}\) in \(J_{1}\). Indeed, for every increasing homeomorphism
\(\lambda\) of \([u_{0},u_{1}]\), the function \(T^{+}_{y}\circ\lambda\) only
takes the values \(1\) and \(2\), whereas, for \(k:=\lceil n/2\rceil\in\{1,\ldots,n-1\}\),
\(T^{+}_{y_{n}}\) takes the value \(1+k/n\in[3/2,3/2+1/n]\). Hence
\[
\sup_{u}\bigl|T^{+}_{y_{n}}(u)-T^{+}_{y}(\lambda(u))\bigr|\geq\frac{1}{2}-\frac{1}{n}
\]
for every such \(\lambda\), so that the \(J_{1}\) distance between
\(T^{+}_{y_{n}}\) and \(T^{+}_{y}\) is at least \(1/2-1/n\). Here the \(J_{1}\)
distance on \(D[u_{0},u_{1}]\) is
\[
d_{J_{1}}(f,h):=\inf_{\lambda}\max\Bigl\{\sup_{u}|\lambda(u)-u|,\ \sup_{u}|f(u)-h(\lambda(u))|\Bigr\},
\]
the infimum being taken over all increasing homeomorphisms \(\lambda\) of
\([u_{0},u_{1}]\) (see~\cite[Section~12]{Billingsley1999}).
\end{example}

\section{Verification routes for \texorpdfstring{\NT}{(NT)} and worked examples}\label{sec:verification}

Sections~\ref{sec:continuity} and~\ref{sec:profiles} establish the
convergence of exit times and of exit-time profiles as continuous-mapping
statements under condition \NT{}, but do not prescribe how \NT{} should be
checked in concrete models. This section records two verification routes
that cover many standard applications, and then presents two worked
examples illustrating the full pipeline: a deterministic moving boundary
(Section~\ref{sec:example-donsker}) and a random domain that is
approximated jointly with the trajectory and coupled to the state
(Section~\ref{sec:example-coupled}).

In this section, time is written either as an argument, as in \(Y(t)\), or,
in statements involving stochastic calculus, as a subscript, as in
\(Y_{t}\); in particular, \(Y_{0}\) denotes an initial value. In the worked examples, the
approximation index of paths, barriers, domains and exit times is written
as a superscript, as in \(X^{(n)}\), \(\Phi^{(n)}\), \(\Gamma^{(n)}\) and
\(\tau^{(n)}\). Theorems~\ref{thm:weak} and~\ref{thm:profiles} are applied
with \(X_{n}:=X^{(n)}\) and \(\Phi_{n}:=\Phi^{(n)}\), and \(\tau^{+}_{n}\)
denotes the profile \eqref{eq:profiles} of \((X^{(n)},\Phi^{(n)})\).

Let \(Y\) be a real càdlàg process, typically a scalarized process
\(Y(t)=\Phi(t,X(t))\), and let \(\tau:=\inf\{t\geq 0:\ Y(t)\geq 0\}\).
Condition \eqref{eq:NTplus} requires that, on \(\{\tau<\infty\}\),
\begin{equation}\label{eq:NTplus-process}
\sup_{\tau\leq s\leq\tau+\varepsilon}Y(s)>0\qquad\text{for every } \varepsilon>0,
\end{equation}
which rules out sticking at the boundary level after first passage.
Condition \eqref{eq:NTminus} holds automatically when \(Y\) is continuous
(Remark~\ref{rem:NTminus}), since a continuous function that stays strictly
below \(0\) on a compact interval attains a strictly negative maximum there.

\subsection{\texorpdfstring{Route A: a jump overshoot implies \eqref{eq:NTplus}}{Route A: a jump overshoot implies (NT+)}}\label{sec:routeA}

This route applies when the scalarized process crosses the boundary by a
jump. If
\begin{equation}\label{eq:overshoot}
Y(\tau)>0\qquad\text{on } \{\tau<\infty\},
\end{equation}
then \eqref{eq:NTplus-process} holds, because right-continuity provides
\(\delta>0\) such that \(Y(t)>Y(\tau)/2>0\) for all \(t\in[\tau,\tau+\delta]\).
This covers many Lévy-driven and compound Poisson settings in which first
passage occurs by overshoot. For processes with jumps, \eqref{eq:NTminus} is
a separate requirement; by Remark~\ref{rem:NTminus}, it holds if and only if
\(Y(t-)<0\) for all \(t\in(0,\tau)\), for instance if \(Y\) has no downward
jumps on \((0,\tau)\).

\subsection{Route B: continuous semimartingales and a non-characteristic Itô criterion}\label{sec:routeB}

The following criterion reduces \eqref{eq:NTplus} to the non-degeneracy of
the martingale part of the scalarized process near the boundary. It is
stated first for continuous semimartingales and then specialized to Itô
processes, for which it becomes a non-characteristic boundary condition.

\begin{proposition}[Continuous semimartingale criterion]\label{prop:semimartingale}
Let \(Y=Y_{0}+M+V\) be a real continuous semimartingale on a filtered
probability space \((\Omega,\mathcal{F},(\mathcal{F}_{t})_{t\geq 0},\PP)\)
satisfying the usual conditions, where \(M\) is a
continuous local martingale and \(V\) is a continuous adapted process of
finite variation, with \(M_{0}=V_{0}=0\). Let \(\tau:=\inf\{t\geq 0:\ Y_{t}\geq 0\}\).
Suppose that, for almost every outcome in \(\{\tau<\infty,\ Y_{\tau}=0\}\),
there exist numbers \(c>0\), \(K\geq 0\) and \(\delta>0\), which may depend on
the outcome, such that
\begin{equation}\label{eq:semimart-cond}
\langle M\rangle_{\tau+t}-\langle M\rangle_{\tau}\geq c^{2}t
\quad\text{and}\quad
|V_{\tau+t}-V_{\tau}|\leq Kt,\qquad t\in[0,\delta].
\end{equation}
Then \(Y\) satisfies \NT{} almost surely: \eqref{eq:NTminus} holds almost
surely, and \eqref{eq:NTplus-process} holds almost surely on
\(\{\tau<\infty\}\).
\end{proposition}

\begin{proof}
Condition \eqref{eq:NTminus} holds because \(Y\) is continuous. On
\(\{\tau<\infty\}\) we have \(Y_{\tau}\geq 0\) by continuity, and on
\(\{\tau<\infty,\ Y_{\tau}>0\}\) (which is contained in \(\{\tau=0\}\)),
\eqref{eq:NTplus-process} holds by continuity. It remains to consider the
event \(\{\tau<\infty,\ Y_{\tau}=0\}\). The process \(N\) defined by
\(N_{t}:=M_{\tau+t}-M_{\tau}\) on \(\{\tau<\infty\}\) and \(N_{t}:=0\) on
\(\{\tau=\infty\}\), \(t\geq 0\), is a continuous local martingale with
respect to \((\mathcal{F}_{\tau+t})_{t\geq 0}\). Indeed, for \(j\in\N\), let
\(\mu_{j}:=\inf\{t\geq 0:\ |M_{t}|\geq j\}\), and let
\(\nu_{j}:=(\mu_{j}-\tau)^{+}\) on \(\{\tau<\infty\}\) and \(\nu_{j}:=\infty\)
on \(\{\tau=\infty\}\). Then \(\nu_{j}\) is an
\((\mathcal{F}_{\tau+t})\)-stopping time, \(\nu_{j}\uparrow\infty\), and
\(N_{t\wedge\nu_{j}}=M^{\mu_{j}}_{\tau+t}-M^{\mu_{j}}_{\tau}\) on
\(\{\tau<\infty\}\), where \(M^{\mu_{j}}_{s}:=M_{s\wedge\mu_{j}}\) denotes the
stopped process; this is a bounded \((\mathcal{F}_{\tau+t})\)-martingale
by the optional stopping theorem applied to the bounded martingale
\(M^{\mu_{j}}\). Moreover, \(\langle N\rangle=\theta\), where
\(\theta(t):=\langle M\rangle_{\tau+t}-\langle M\rangle_{\tau}\) on
\(\{\tau<\infty\}\) and \(\theta(t):=0\) on \(\{\tau=\infty\}\).
By the Dambis--Dubins--Schwarz theorem~\cite[Section~V.1]{RevuzYor1999}, on
an enlargement of the probability space there is a standard Brownian motion
\(\widehat{B}\) such that \(N_{t}=\widehat{B}_{\theta(t)}\) for all \(t\geq 0\).
By the law of the iterated logarithm~\cite[Section~II.1]{RevuzYor1999},
\(\limsup_{s\downarrow 0}\widehat{B}_{s}/s=+\infty\) almost surely.

Now fix an outcome in the event \(\{\tau<\infty,\ Y_{\tau}=0\}\) such that
\eqref{eq:semimart-cond} holds and
\(\limsup_{s\downarrow 0}\widehat{B}_{s}/s=+\infty\), and let
\(\varepsilon\in(0,\delta]\). The function \(\theta\) is continuous and
nondecreasing, with \(\theta(0)=0\) and \(\theta(\varepsilon)\geq c^{2}\varepsilon>0\).
Choose \(s\in(0,\theta(\varepsilon)]\) with \(\widehat{B}_{s}>(K/c^{2})s\); by
the intermediate value theorem there is \(t\in(0,\varepsilon]\) with
\(\theta(t)=s\), and \eqref{eq:semimart-cond} gives \(c^{2}t\leq\theta(t)=s\).
Hence
\[
Y_{\tau+t}=Y_{\tau+t}-Y_{\tau}=\widehat{B}_{s}+V_{\tau+t}-V_{\tau}\geq\widehat{B}_{s}-Kt\geq\widehat{B}_{s}-\frac{K}{c^{2}}\,s>0 .
\]
As \(\varepsilon\in(0,\delta]\) is arbitrary, \eqref{eq:NTplus-process}
follows.
\end{proof}

\begin{corollary}[Itô criterion: non-characteristic boundary]\label{cor:ito}
Let \(k,m\geq 1\), let \(W\) be an \(m\)-dimensional Brownian motion, and let
\(X\) be a continuous \(\R^{k}\)-valued Itô process such that
\begin{equation}\label{eq:ito-sde}
dX_{t}=b(t,X_{t})\,dt+\sigma(t,X_{t})\,dW_{t},
\end{equation}
where \(b:[0,\infty)\times\R^{k}\to\R^{k}\) and
\(\sigma:[0,\infty)\times\R^{k}\to\R^{k\times m}\) are locally bounded Borel
functions (\(\sigma\) may be degenerate). Let
\(\Phi\in C^{1,2}([0,\infty)\times\R^{k})\) be a deterministic barrier, let
\(Y_{t}:=\Phi(t,X_{t})\), and let \(\tau:=\inf\{t\geq 0:\ Y_{t}\geq 0\}\).
Suppose that, for all \(T,R>0\), there exist \(\eta>0\) and \(c>0\) such that
the non-characteristic condition
\begin{equation}\label{eq:noncharacteristic}
\bigl\|\sigma(t,z)^{\top}\nabla_{z}\Phi(t,z)\bigr\|\geq c
\end{equation}
holds on the space--time neighborhood
\begin{equation}\label{eq:NTeta}
\cN_{T,R,\eta}:=\{(t,z)\in[0,T]\times\mathbb{B}_{R}:\ |\Phi(t,z)|\leq\eta\}
\end{equation}
of the boundary, where, in this corollary, \(\mathbb{B}_{R}\) denotes the closed ball of
radius \(R\) in \(\R^{k}\), and \NT{} is understood with \(d\) replaced by \(k\). Then \((X,\Phi)\) satisfies \NT{} almost surely. The same
conclusion holds at a level \(u\in\R\) if the assumption holds with \(\Phi\)
replaced by \(\Phi-u\).
\end{corollary}

\begin{proof}
By Itô's formula, \(Y=Y_{0}+M+V\) with
\(M_{t}:=\int_{0}^{t}\beta_{s}^{\top}\,dW_{s}\) and
\(V_{t}:=\int_{0}^{t}\alpha_{s}\,ds\), where
\(\beta_{s}:=\sigma(s,X_{s})^{\top}\nabla_{z}\Phi(s,X_{s})\) and
\[
\alpha_{s}:=\Bigl(\partial_{t}\Phi+b\cdot\nabla_{z}\Phi+\tfrac{1}{2}\operatorname{tr}\bigl(\sigma\sigma^{\top}\nabla_{z}^{2}\Phi\bigr)\Bigr)(s,X_{s}).
\]
In particular, \(\langle M\rangle_{t}=\int_{0}^{t}\|\beta_{s}\|^{2}\,ds\). Fix an
outcome in \(\{\tau<\infty,\ Y_{\tau}=0\}\). Let \(T>\tau+1\) and
\(R\geq\sup_{r\leq T}\|X_{r}\|\) be integers, and let \(\eta\) and \(c\) be the
constants provided by the assumption for this pair \((T,R)\). By continuity
of \(Y\), there is \(\delta\in(0,1]\) such that \(|Y_{\tau+t}|\leq\eta\), hence
\((\tau+t,X_{\tau+t})\in\cN_{T,R,\eta}\) and \(\|\beta_{\tau+t}\|\geq c\), for all
\(t\in[0,\delta]\); therefore
\(\langle M\rangle_{\tau+t}-\langle M\rangle_{\tau}\geq c^{2}t\) on \([0,\delta]\).
Moreover, \((s,X_{s})\in[0,T]\times\mathbb{B}_{R}\) for \(s\leq T\), a compact
set on which \(b\) and \(\sigma\) are bounded and the derivatives of \(\Phi\)
are continuous, hence bounded. Thus \(|\alpha_{s}|\leq K\) for \(s\leq T\)
and some finite \(K\), and \(|V_{\tau+t}-V_{\tau}|\leq Kt\) for
\(t\in[0,\delta]\). The conditions of Proposition~\ref{prop:semimartingale}
are therefore satisfied almost surely on \(\{\tau<\infty,\ Y_{\tau}=0\}\),
and \NT{} follows. The last statement follows by applying the result to
\(\Phi-u\).
\end{proof}

\subsection{Worked example I: Donsker approximation and a deterministic moving boundary}\label{sec:example-donsker}

This example illustrates the complete recipe: joint convergence of
paths and barriers, verification of \NT{} for the limit, and conclusions
for exit times and exit-time profiles via Theorems~\ref{thm:weak}
and~\ref{thm:profiles}.

\paragraph{Setup}
Let \((\xi_{k})_{k\geq 1}\) be i.i.d.\ real random variables with
\(\E[\xi_{1}]=0\) and \(\E[\xi_{1}^{2}]=1\), let \(S_{n}:=\sum_{k=1}^{n}\xi_{k}\),
and let
\begin{equation}\label{eq:donsker-walk}
X^{(n)}(t):=\frac{1}{\sqrt{n}}\,S_{\lfloor nt\rfloor},\qquad t\geq 0.
\end{equation}
By Donsker's invariance principle, \(X^{(n)}\Rightarrow B\) in \(\Done\), where
\(B\) is a standard Brownian motion (see, e.g.,~\cite[Sections~14 and~16]{Billingsley1999} or~\cite[Chapter~VII]{JacodShiryaev2003}). Let
\(g\in C^{1}([0,\infty))\) be a deterministic moving boundary, consider the
time-dependent domain \(\Gamma_{t}:=(-\infty,g(t))\), and represent it by the
barrier
\[
\Phi(t,z):=z-g(t),\qquad\text{so that}\qquad\Gamma_{t}=\{z\in\R:\ \Phi(t,z)<0\}.
\]
The exit times are
\[
\tau^{(n)}:=\inf\{t\geq 0:\ X^{(n)}(t)\geq g(t)\},\qquad
\tau:=\inf\{t\geq 0:\ B(t)\geq g(t)\}.
\]

\paragraph{Verification of \texorpdfstring{\NT}{(NT)} for the limit}
The scalarized limit process \(Y_{t}:=\Phi(t,B_{t})\) is given by
\(Y_{t}=B_{t}-g(t)\), so that \(dY_{t}=-g'(t)\,dt+dB_{t}\). Corollary~\ref{cor:ito} applies with \(k=m=1\),
\(b\equiv 0\), \(\sigma\equiv 1\) and \(\nabla_{z}\Phi\equiv 1\), so that
\eqref{eq:noncharacteristic} holds with \(c=1\). Hence \((B,\Phi)\)
satisfies \NT{} almost surely. Applying Corollary~\ref{cor:ito} to \(\Phi-u\)
in the same way shows that, for each fixed level \(u\in\R\), \((B,\Phi)\)
satisfies \NT{} at level \(u\) almost surely.

\paragraph{Exit-time convergence}
Take \(\Phi^{(n)}:=\Phi\) for all \(n\). Since the second component is
deterministic, \((X^{(n)},\Phi)\Rightarrow(B,\Phi)\) jointly~\cite[Theorem~3.9]{Billingsley1999}, and Theorem~\ref{thm:weak} yields
\(\tau^{(n)}\Rightarrow\tau\) in \([0,\infty]\).

\paragraph{\texorpdfstring{Profile convergence in \(M_{1}\)}{Profile convergence in M1}}
Fix \(u_{0}<u_{1}\). For \(u\in[u_{0},u_{1}]\), the profiles \eqref{eq:profiles}
of \((X^{(n)},\Phi)\) and \((B,\Phi)\) are
\[
\begin{aligned}
\tau^{+}_{n}(u)&=\inf\{t\geq 0:\ X^{(n)}(t)>g(t)+u\},\\
\tau^{+}(u)&=\inf\{t\geq 0:\ B(t)>g(t)+u\}.
\end{aligned}
\]
Since \NT{} holds almost surely at each of the two fixed levels \(u_{0}\)
and \(u_{1}\), hence almost surely at both, Theorem~\ref{thm:profiles} yields \(\tau^{+}_{n}\Rightarrow\tau^{+}\) in
\(\Ms[u_{0},u_{1}]\). If \(g\) is bounded above, then \(\tau^{+}(u_{1})<\infty\)
almost surely, because \(\limsup_{t\to\infty}B_{t}=+\infty\), and the
convergence also holds in \((D[u_{0},u_{1}],M_{1})\) in the sense of
Remark~\ref{rem:finite-profiles}.

\subsection{Worked example II: a random domain coupled to the state}\label{sec:example-coupled}

In Section~\ref{sec:example-donsker}, the domain is deterministic and is not
approximated: only the trajectory is perturbed. We now treat the situation
announced in the abstract and in Remark~\ref{rem:dependence}. The
approximating domains are random; they are perturbed jointly with the
trajectory, because they are built from the same increments; and they are
coupled to the state. The limiting domain is itself random and depends on
the limiting trajectory.

\subsubsection{Setup: a tracking domain centered at a moving average}\label{sec:example-coupled-setup}

Fix \(d\geq 1\) and let \((\xi_{k})_{k\geq 1}\) be i.i.d.\ \(\R^{d}\)-valued
random vectors with \(\E[\xi_{1}]=0\) and \(\E[\xi_{1}\xi_{1}^{\top}]=I_{d}\).
As in \eqref{eq:donsker-walk}, let
\begin{equation}\label{eq:walk-d}
X^{(n)}(t):=\frac{1}{\sqrt{n}}\sum_{k=1}^{\lfloor nt\rfloor}\xi_{k},\qquad t\geq 0.
\end{equation}
By the multidimensional version of Donsker's invariance principle (see,
e.g., \cite[Chapter~VII]{JacodShiryaev2003}), \(X^{(n)}\Rightarrow W\) in
\(\Dd\), where \(W\) is a standard \(d\)-dimensional Brownian motion. The
domains are built from three ingredients.

\emph{(a) Moving average.} Fix \(\kappa>0\). For \(x\in\Dd\), let
\begin{equation}\label{eq:ewma}
\cE(x)(t):=\kappa\int_{0}^{t}e^{-\kappa(t-s)}\,x(s)\,ds,\qquad t\geq 0,
\end{equation}
be the exponentially weighted moving average of \(x\). The function
\(\cE(x)\) is continuous, \(\cE(x)(0)=0\), and
\(\cE(x)(t)=\kappa\int_{0}^{t}\bigl(x(s)-\cE(x)(s)\bigr)\,ds\): the average
relaxes toward the path at rate \(\kappa\).

\emph{(b) Tolerance.} Fix \(\rho\in C^{1,2}([0,\infty)\times\R^{d})\) with
\(0<\rho_{*}\leq\rho\leq\rho^{*}<\infty\) for some constants \(\rho_{*}\) and
\(\rho^{*}\). The admissible deviation at time \(t\), when the moving average
is at \(a\), is \(\rho(t,a)\); it may depend both on time and on the moving
average itself.

\emph{(c) Calibration.} Fix \(\Tcal>0\) and, for \(n\geq 1/\Tcal\), let
\(N_{n}:=\lfloor n\Tcal\rfloor\) and
\begin{equation}\label{eq:calibration}
\hat{s}_{n}:=\Bigl(\frac{1}{d\,N_{n}}\sum_{k=1}^{N_{n}}\|\xi_{k}\|^{2}\Bigr)^{1/2}.
\end{equation}
Since \(\xi_{k}=\sqrt{n}\,\Delta X^{(n)}(k/n)\), where \(\Delta\) denotes the
jump (Section~\ref{sec:topologies}), we have
\(\hat{s}_{n}^{2}=(n/(d\,N_{n}))\sum_{0<s\leq\Tcal}\|\Delta X^{(n)}(s)\|^{2}\):
the quantity \(\hat{s}_{n}\) is a normalized realized volatility of the path
\(X^{(n)}\) on the calibration window \([0,\Tcal]\), and \(\hat{s}_{n}\to 1\)
almost surely (Step~1 of the proof of Proposition~\ref{prop:coupled} below).

Let \(A^{(n)}:=\cE(X^{(n)})\) and \(A:=\cE(W)\). The approximating and
limiting domains are the moving balls
\begin{align}
\Gamma^{(n)}_{t}&:=\bigl\{z\in\R^{d}:\ \|z-A^{(n)}(t)\|<\hat{s}_{n}\,\rho\bigl(t,A^{(n)}(t)\bigr)\bigr\},\label{eq:domain-n}\\
\Gamma_{t}&:=\bigl\{z\in\R^{d}:\ \|z-A(t)\|<\rho\bigl(t,A(t)\bigr)\bigr\},\label{eq:domain-limit}
\end{align}
with barrier representations
\begin{equation}\label{eq:coupled-barriers}
\begin{aligned}
\Phi^{(n)}(t,z)&:=\|z-A^{(n)}(t)\|^{2}-\hat{s}_{n}^{2}\,\rho\bigl(t,A^{(n)}(t)\bigr)^{2},\\
\Phi(t,z)&:=\|z-A(t)\|^{2}-\rho\bigl(t,A(t)\bigr)^{2}.
\end{aligned}
\end{equation}
Each \(\Phi^{(n)}\) is a random element of \(\cH\): although \(X^{(n)}\) is a
step function, \(A^{(n)}\) is continuous, and so is \(\Phi^{(n)}\)
(measurability follows from Step~1 of the proof of
Proposition~\ref{prop:coupled}). The corresponding exit times
\(\tau^{(n)}:=\tau(X^{(n)},\Phi^{(n)})\) and \(\tau:=\tau(W,\Phi)\) are
\begin{equation}\label{eq:coupled-exit}
\begin{aligned}
\tau^{(n)}&=\inf\bigl\{t\geq 0:\ \|X^{(n)}(t)-A^{(n)}(t)\|\geq\hat{s}_{n}\,\rho\bigl(t,A^{(n)}(t)\bigr)\bigr\},\\
\tau&=\inf\bigl\{t\geq 0:\ \|W(t)-A(t)\|\geq\rho\bigl(t,A(t)\bigr)\bigr\}.
\end{aligned}
\end{equation}
Thus \(\Gamma^{(n)}\) is a tracking domain: its center follows the state
with an exponential lag, its size is prescribed by \(\rho\) and calibrated
on the observed data, and \(\tau^{(n)}\) is the first time the deviation
of the state from its own moving average reaches or exceeds the calibrated
tolerance. For \(d=1\) and constant \(\rho\), \(\Gamma^{(n)}\) is a band
around the moving average.

In contrast with Section~\ref{sec:example-donsker}, the domain is now
coupled to the state in two different ways. First, the center
\(A^{(n)}(t)\), which also enters the tolerance \(\rho(t,A^{(n)}(t))\), is an
adapted and continuous functional of the path. Second, the factor
\(\hat{s}_{n}\) depends on the increments of \(X^{(n)}\) up to time
\(N_{n}/n\), so that \(\Gamma^{(n)}_{t}\) need not be adapted to the natural
filtration of \(X^{(n)}\) for \(t<N_{n}/n\). The domain \(\Gamma^{(n)}\) is not
adapted when
\(\|\xi_{1}\|^{2}\) is not almost surely constant, for instance for the
centered exponential increments of Remark~\ref{rem:coupled}(v): since
\(\rho>0\) and \(A^{(n)}\) is adapted, adaptedness at a time \(t<N_{n}/n\)
would make \(\hat{s}_{n}^{2}\), and hence \(\|\xi_{N_{n}}\|^{2}\), measurable
with respect to the \(\sigma\)-field generated by
\(\xi_{1},\dots,\xi_{N_{n}-1}\), of which
\(\|\xi_{N_{n}}\|^{2}\) is independent, so that \(\|\xi_{N_{n}}\|^{2}\) would
be almost surely constant. By contrast, if the coordinates of \(\xi_{1}\)
are independent Rademacher variables, then \(\hat{s}_{n}=1\) and the domains
are adapted. Moreover, \(\hat{s}_{n}\) is obtained from \(X^{(n)}\) through a
functional that is not continuous for the \(J_{1}\) topology
(Remark~\ref{rem:coupled}(ii)). Finally, the limiting
domain \(\Gamma\) is random and is a functional of \(W\).

\subsubsection{Joint convergence, verification of \texorpdfstring{\NT}{(NT)}, and conclusions}\label{sec:example-coupled-results}

\begin{proposition}[Exit from a random tracking domain]\label{prop:coupled}
In the setting \eqref{eq:walk-d}--\eqref{eq:coupled-exit}, the following hold.
\begin{enumerate}[(a)]
\item \((X^{(n)},\Phi^{(n)})\Rightarrow(W,\Phi)\) in \(\Dd\times\cH\).
\item Almost surely, \(\tau<\infty\) and \((W,\Phi)\) satisfies \NT{}.
\item \(\tau^{(n)}\Rightarrow\tau\) in \([0,\infty]\).
\item Let \(u_{0}<u_{1}\) with \(u_{0}>-\rho_{*}^{2}\), and let
\(\tau^{+}_{n}\) and \(\tau^{+}\) be the exit-time profiles
\eqref{eq:profiles} of \((X^{(n)},\Phi^{(n)})\) and \((W,\Phi)\) on
\([u_{0},u_{1}]\). Then \(\tau^{+}_{n}\Rightarrow\tau^{+}\) in
\(\Ms[u_{0},u_{1}]\). Moreover, \(\tau^{+}(u_{1})<\infty\) almost surely, so
that the convergence also holds in \((D[u_{0},u_{1}],M_{1})\) in the sense of
Remark~\ref{rem:finite-profiles}.
\item If \(\rho\equiv r_{0}\) is constant, then, for every \(\lambda>0\), with the convention \(e^{-\lambda\cdot\infty}:=0\),
\begin{equation}\label{eq:kummer-limit}
\lim_{n\to\infty}\E\bigl[e^{-\lambda\tau^{(n)}}\bigr]=\E\bigl[e^{-\lambda\tau}\bigr]
=\frac{1}{\hyp\bigl(\lambda/(2\kappa);d/2;\kappa r_{0}^{2}\bigr)},
\end{equation}
where
\[
\hyp(\alpha;\gamma;\zeta):=\sum_{j=0}^{\infty}\frac{(\alpha)_{j}}{(\gamma)_{j}\,j!}\,\zeta^{j},
\]
with \((\alpha)_{0}:=1\) and \((\alpha)_{j}:=\alpha(\alpha+1)\cdots(\alpha+j-1)\),
is Kummer's confluent hypergeometric function, denoted by
\(M(\alpha,\gamma,\zeta)\) in~\cite[(13.2.2)]{Olver2010}. In particular,
\(\E[e^{-\lambda\tau}]=e^{-\kappa r_{0}^{2}}\) for \(\lambda:=d\kappa\).
\end{enumerate}
\end{proposition}

\begin{proof}
\emph{Step 1: joint convergence.} For \(x,x'\in\Dd\) and \(T>0\),
\begin{equation}\label{eq:ewma-lipschitz}
\sup_{t\leq T}\bigl\|\cE(x)(t)-\cE(x')(t)\bigr\|\leq\kappa\int_{0}^{T}\|x(s)-x'(s)\|\,ds .
\end{equation}
If \(x_{m}\to x\) in \(\Dd\), then \(x_{m}(s)\to x(s)\) at every continuity
point \(s\) of \(x\), hence for Lebesgue-almost every \(s\), and
\(\sup_{m}\sup_{s\leq T}\|x_{m}(s)\|<\infty\) for every \(T\)
(Section~\ref{sec:topologies}). By dominated convergence and
\eqref{eq:ewma-lipschitz}, \(\cE(x_{m})\to\cE(x)\) uniformly on compacts, so
that \(\cE\) is continuous from \(\Dd\) to \(C([0,\infty),\R^{d})\). Next, define
\(G:\Dd\times\R\to\Dd\times\cH\) by \(G(x,q):=(x,\Phi^{x,q})\), where
\[
\Phi^{x,q}(t,z):=\|z-\cE(x)(t)\|^{2}-q^{2}\rho\bigl(t,\cE(x)(t)\bigr)^{2}.
\]
Let \((x_{m},q_{m})\to(x,q)\), fix \(T,R>0\), write \(a_{m}:=\cE(x_{m})\) and
\(a:=\cE(x)\), and let \(K\subset\R^{d}\) be a compact set containing
\(a_{m}(t)\) and \(a(t)\) for all \(m\) and all \(t\leq T\). Using the identity
\(\|z-v'\|^{2}-\|z-v\|^{2}=(v-v')\cdot(2z-v-v')\) for \(v,v'\in\R^{d}\), we
get, for \(t\leq T\) and \(\|z\|\leq R\),
\[
\begin{aligned}
\bigl|\Phi^{x_{m},q_{m}}(t,z)-\Phi^{x,q}(t,z)\bigr|
\leq{}&\|a_{m}(t)-a(t)\|\,\Bigl(2R+2\max_{w\in K}\|w\|\Bigr)\\
&+|q_{m}^{2}-q^{2}|\,(\rho^{*})^{2}\\
&+q^{2}\,\bigl|\rho\bigl(t,a_{m}(t)\bigr)^{2}-\rho\bigl(t,a(t)\bigr)^{2}\bigr|.
\end{aligned}
\]
The right-hand side tends to \(0\) uniformly in \((t,z)\in[0,T]\times\mathbb{B}_{R}\),
because \(a_{m}\to a\) uniformly on \([0,T]\) and \(\rho^{2}\) is uniformly
continuous on \([0,T]\times K\). Hence \(G\) is continuous, and in
particular Borel. By the strong law of large numbers,
\(\hat{s}_{n}^{2}\to\E[\|\xi_{1}\|^{2}]/d=1\) almost surely. Since the limit is
a constant and \(X^{(n)}\Rightarrow W\), we obtain
\((X^{(n)},\hat{s}_{n})\Rightarrow(W,1)\) in \(\Dd\times\R\)~\cite[Theorem~3.9]{Billingsley1999}. As
\((X^{(n)},\Phi^{(n)})=G(X^{(n)},\hat{s}_{n})\) and \((W,\Phi)=G(W,1)\), the
continuous mapping theorem~\cite[Theorem~2.7]{Billingsley1999} gives (a).

\emph{Step 2: \NT{} by augmentation of the state.} The limiting barrier
\(\Phi\) is random, so that Corollary~\ref{cor:ito}, which is formulated
for deterministic barriers, does not apply to \((W,\Phi)\) directly.
However, the randomness of \(\Phi\) is carried by the finite-dimensional
process \(A\), which suggests augmenting the state. Let
\(A_{t}:=A(t)\) and \(Z_{t}:=(W_{t},A_{t})\in\R^{2d}\). Since
\(A_{t}=\kappa\int_{0}^{t}(W_{s}-A_{s})\,ds\), \(Z\) is a continuous Itô process
of the form \eqref{eq:ito-sde} with \(k=2d\), \(m=d\), driving Brownian
motion \(W\) and constant diffusion matrix \(\sigma\in\R^{2d\times d}\):
\begin{equation}\label{eq:augmented}
dZ_{t}=b(Z_{t})\,dt+\sigma\,dW_{t},\qquad
b(z,a):=\begin{pmatrix}0\\ \kappa(z-a)\end{pmatrix},\qquad
\sigma:=\begin{pmatrix}I_{d}\\ 0\end{pmatrix}.
\end{equation}
On the augmented space the barrier becomes deterministic:
\begin{equation}\label{eq:augmented-barrier}
\Phi(t,W_{t})=\Psi(t,Z_{t}),\qquad\Psi(t,z,a):=\|z-a\|^{2}-\rho(t,a)^{2},
\end{equation}
where \(\Psi\in C^{1,2}([0,\infty)\times\R^{2d})\) is deterministic. Hence the
scalarized limit process \(Y_{t}:=\Phi(t,W_{t})=\Psi(t,Z_{t})\) and the exit
time \(\tau=\tau(W,\Phi)=\tau(Z,\Psi)\) are the same in both formulations.
Since \NT{} only involves the scalarized path, it suffices to verify it for
\((Z,\Psi)\). Now
\(\sigma^{\top}\nabla_{(z,a)}\Psi(t,z,a)=\nabla_{z}\Psi(t,z,a)=2(z-a)\).
Let \(\eta:=\rho_{*}^{2}/2\). On the set
\(\{(t,z,a):\ |\Psi(t,z,a)|\leq\eta\}\), which contains the sets
\(\cN_{T,R,\eta}\) of \eqref{eq:NTeta} (with \(\Psi\) in place of \(\Phi\) and \(k=2d\)) for all \(T,R>0\), we have
\(\|z-a\|^{2}\geq\rho(t,a)^{2}-\eta\geq\rho_{*}^{2}/2\), so that
\begin{equation}\label{eq:augmented-noncharacteristic}
\bigl\|\sigma^{\top}\nabla_{(z,a)}\Psi(t,z,a)\bigr\|=2\|z-a\|\geq\sqrt{2}\,\rho_{*}>0 .
\end{equation}
Since \(b\) and \(\sigma\) are locally bounded, Corollary~\ref{cor:ito} shows
that \((Z,\Psi)\), and hence \((W,\Phi)\), satisfies \NT{} almost surely. Note
that the augmented diffusion \eqref{eq:augmented} is degenerate, as
\(\sigma\sigma^{\top}\) has rank \(d<2d\). This is no obstacle: the
non-characteristic condition \eqref{eq:noncharacteristic} only requires the
noise to have a nonvanishing component along the normal to the level sets
of the barrier, and the \(z\)-component \(2(z-a)\) of \(\nabla_{(z,a)}\Psi\)
does not vanish near the boundary.

It remains to show that \(\tau<\infty\) almost surely. The deviation
\(U:=W-A\) is an Ornstein--Uhlenbeck process,
\(dU_{t}=dW_{t}-\kappa U_{t}\,dt\), \(U_{0}=0\), and
\({Y_{t}=\|U_{t}\|^{2}-\rho(t,A_{t})^{2}}\). For every \(u\in\R\),
\(\|U_{t}\|^{2}\geq(\rho^{*})^{2}+|u|\) implies \(Y_{t}\geq u\), so it suffices
to show that \(\sup_{t\geq 0}\|U_{t}\|=\infty\) almost surely. The first
coordinate \(U^{1}\) of \(U\) satisfies
\(e^{\kappa t}U_{t}^{1}=\int_{0}^{t}e^{\kappa s}\,dW_{s}^{1}\), a continuous
martingale with quadratic variation
\(\vartheta(t):=(e^{2\kappa t}-1)/(2\kappa)\to\infty\). By the
Dambis--Dubins--Schwarz theorem~\cite[Section~V.1]{RevuzYor1999}, there is a
one-dimensional standard Brownian motion \(\widetilde{W}\) such that
\(\int_{0}^{t}e^{\kappa s}\,dW_{s}^{1}=\widetilde{W}_{\vartheta(t)}\), whence
\(U_{t}^{1}=\widetilde{W}_{\vartheta(t)}/(1+2\kappa\vartheta(t))^{1/2}\). The law of
the iterated logarithm for \(\widetilde{W}\)~\cite[Section~II.1]{RevuzYor1999}
gives \({\sup_{t\geq 0}U_{t}^{1}=+\infty}\) almost surely. This proves (b) and
shows that, almost surely, \({\sup_{t\geq 0}Y_{t}=+\infty}\), so that
\(T^{+}_{Y}(u)<\infty\), and a fortiori \(\tau_{u}(W,\Phi)<\infty\), for all
\(u\in\R\).

\emph{Step 3: exit times.} By (a), (b) and Theorem~\ref{thm:weak},
\(\tau^{(n)}\Rightarrow\tau\) in \([0,\infty]\), which is (c).

\emph{Step 4: exit-time profiles.} Fix a deterministic level
\(u>-\rho_{*}^{2}\). We show that \((W,\Phi)\) satisfies \NT{} at level \(u\)
almost surely. Let \(\eta_{u}:=(\rho_{*}^{2}+u)/2>0\). On the set
\(\{(t,z,a):\ |\Psi(t,z,a)-u|\leq\eta_{u}\}\) we have
\(\|z-a\|^{2}\geq\rho_{*}^{2}+u-\eta_{u}=\eta_{u}\), hence
\(\|\sigma^{\top}\nabla_{(z,a)}\Psi\|\geq 2\eta_{u}^{1/2}>0\), and
Corollary~\ref{cor:ito}, applied to \((Z,\Psi-u)\), yields the claim.
Applying this to the two fixed levels \(u_{0}\) and \(u_{1}\) and
intersecting the two resulting events of probability one, we obtain that,
almost surely, \((W,\Phi)\) satisfies \NT{} at the levels \(u_{0}\) and
\(u_{1}\). Together with (a), this verifies the hypotheses of
Theorem~\ref{thm:profiles}, which gives the first statement in (d); the second follows from the last sentence of
Step~2. We do not claim that, almost surely, \NT{} holds simultaneously at
all levels \(u>-\rho_{*}^{2}\); in general it does not
(Remark~\ref{rem:end-levels}). The condition
\(u>-\rho_{*}^{2}\) ensures that the level-\(u\) domains are balls of radius
at least \((\rho_{*}^{2}+u)^{1/2}>0\), which keeps the level set
\(\{\Psi=u\}\) away from the set \(\{z=a\}\) on which
\(\sigma^{\top}\nabla_{(z,a)}\Psi\) vanishes; it is sufficient, not necessary.

\emph{Step 5: explicit limit law.} Let \(\rho\equiv r_{0}\), so that
\(\tau=\inf\{t\geq 0:\ \|U_{t}\|\geq r_{0}\}\) is the exit time of the
Ornstein--Uhlenbeck process \(U\) from the ball of radius \(r_{0}\) centered
at \(0\), i.e.\ the first hitting time of \(r_{0}^{2}\) by the
Cox--Ingersoll--Ross-type process \(\|U\|^{2}\). Fix \(\lambda>0\) and let
\(\varphi(\zeta):=\hyp(\lambda/(2\kappa);d/2;\zeta)\). This entire function
solves Kummer's equation
\(\zeta\varphi''+(d/2-\zeta)\varphi'-(\lambda/(2\kappa))\varphi=0\)~\cite[(13.2.1)]{Olver2010}, and a direct computation shows that
\(f(z):=\varphi(\kappa\|z\|^{2})\) satisfies
\(\frac{1}{2}\operatorname{tr}\nabla^{2}f(z)-\kappa\,z\cdot\nabla f(z)=\lambda f(z)\) on \(\R^{d}\).
By Itô's formula, \(e^{-\lambda(t\wedge\tau)}f(U_{t\wedge\tau})\) is a
martingale, bounded by \(\max_{\|z\|\leq r_{0}}|f(z)|\). Since \(\tau<\infty\)
almost surely, optional stopping and bounded convergence give
\(\varphi(\kappa r_{0}^{2})\,\E[e^{-\lambda\tau}]=f(0)=1\), which is the second
equality in \eqref{eq:kummer-limit}; formulas of this type are classical
for hitting times of Cox--Ingersoll--Ross processes (see, e.g.,~\cite{Linetsky2004}). As \(t\mapsto e^{-\lambda t}\), with
\(e^{-\lambda\cdot\infty}:=0\), is bounded and continuous on \([0,\infty]\), the
first equality follows from (c). Finally,
\(\hyp(d/2;d/2;\zeta)=e^{\zeta}\)~\cite[(13.6.1)]{Olver2010}.
\end{proof}

\subsubsection{Comments and a variant with a time-discontinuous coupling}\label{sec:example-coupled-comments}

\begin{remark}[What the example illustrates]\label{rem:coupled}
We comment on several features of Proposition~\ref{prop:coupled}.
\begin{enumerate}[(i)]
\item Neither independence between trajectory and domain nor adaptedness
of the domain to the natural filtration of the trajectory is needed. The domains \(\Gamma^{(n)}\) are built from the same increments as
\(X^{(n)}\) and may anticipate them through \(\hat{s}_{n}\): they do so
whenever \(\|\xi_{1}\|^{2}\) is not almost surely constant
(Section~\ref{sec:example-coupled-setup}), and then, in general,
\(\tau^{(n)}\) is not a stopping time of the natural filtration of
\(X^{(n)}\). The limiting domain \(\Gamma\) is a functional of \(W\). Only
the joint convergence in Proposition~\ref{prop:coupled}(a) and \NT{} for the
limit are used, in line with Remark~\ref{rem:dependence}.
\item Joint convergence need not follow from the continuous mapping theorem
applied to the path alone. The functional
\(x\mapsto\sum_{0<s\leq\Tcal}\|\Delta x(s)\|^{2}\), from which \(\hat{s}_{n}\) is
computed, is not continuous for the \(J_{1}\) topology. Indeed, its value at
\(X^{(n)}\) is \(d\,N_{n}\hat{s}_{n}^{2}/n\to d\,\Tcal>0\), whereas it vanishes at
every continuous path; if it were continuous at almost every path of \(W\),
the continuous mapping theorem would force its values at \(X^{(n)}\) to
converge in distribution to \(0\). Joint convergence is obtained from
Slutsky's lemma instead~\cite[Theorem~3.9]{Billingsley1999}.
\item When the coupling between a random barrier and the state is carried by
a finite-dimensional process, \NT{} can be verified with the
deterministic-barrier criteria of Section~\ref{sec:routeB} after augmenting
the state with that process, provided the augmented state is a continuous
semimartingale whose scalarized path satisfies the assumptions of
Proposition~\ref{prop:semimartingale} (for instance, an Itô process of the
form \eqref{eq:ito-sde}, via Corollary~\ref{cor:ito}). The augmented
diffusion is typically degenerate, and it is the non-characteristic
condition \eqref{eq:noncharacteristic}, not ellipticity, that is required.
The augmentation could also be performed before the limit: one has
\(\tau^{(n)}=\tau\bigl((X^{(n)},A^{(n)}),\Psi^{(n)}\bigr)\) with
\(\Psi^{(n)}(t,z,a):=\|z-a\|^{2}-\hat{s}_{n}^{2}\rho(t,a)^{2}\), whose only
remaining randomness is the scalar \(\hat{s}_{n}\), and
\((X^{(n)},A^{(n)})\Rightarrow(W,A)\) in \(\mathbb{D}^{2d}\) because
\(x\mapsto(x,\cE(x))\) is continuous from \(\Dd\) to \(\mathbb{D}^{2d}\). We
have kept the formulation in \(\R^{d}\), where the domains actually live, to
show that Theorems~\ref{thm:weak} and~\ref{thm:profiles} apply directly to
random, state-coupled and possibly nonadapted barriers; augmentation is only
used as a tool to verify \NT{} for the limit.
\item Random tolerance fields independent of the walk can be handled with
the same proof. Suppose that \(\rho\) is replaced by random tolerance fields
\(\tilde{\rho}^{(n)}\), independent of
\((\xi_{k})_{k\geq 1}\) and nonnegative (that is, almost surely,
\(\tilde{\rho}^{(n)}\geq 0\) on \([0,\infty)\times\R^{d}\)), such that \(\tilde{\rho}^{(n)}\Rightarrow\tilde{\rho}\)
in \(\cH\), where \(\tilde{\rho}\) takes values in \(C^{1,2}([0,\infty)\times\R^{d})\) and
satisfies, almost surely, \(\rho_{*}\leq\tilde{\rho}\leq\rho^{*}\) on
\([0,\infty)\times\R^{d}\) for deterministic constants
\(0<\rho_{*}\leq\rho^{*}<\infty\). Nonnegativity ensures that the domains
\eqref{eq:domain-n}, with \(\rho\) replaced by \(\tilde{\rho}^{(n)}\), are still
represented by the barriers \eqref{eq:coupled-barriers}, because
\(\|z-a\|<r\) if and only if \(\|z-a\|^{2}<r^{2}\) when \(r\geq 0\); a
negative tolerance would give an empty domain, but a nonempty barrier
domain when \(\hat{s}_{n}>0\). Then \((X^{(n)},\hat{s}_{n},\tilde{\rho}^{(n)})\Rightarrow(W,1,\tilde{\rho})\),
with \(W\) and \(\tilde{\rho}\) independent, by~\cite[Theorems~2.8 and~3.9]{Billingsley1999}.
The map \(G\) of Step~1 is also continuous in the tolerance field, and
Steps~2--4 apply conditionally on \(\tilde{\rho}\), because \(W\) is a Brownian motion
independent of \(\tilde{\rho}\). Hence Proposition~\ref{prop:coupled}(a)--(d) remain
valid with \(\rho\) replaced by \(\tilde{\rho}^{(n)}\) and \(\tilde{\rho}\).
\item When \(\rho\equiv r_{0}\), the distance \(\|X^{(n)}-A^{(n)}\|\) is
nonincreasing between the jump times of \(X^{(n)}\), so that \(\tau^{(n)}\)
takes values in the grid \(\{k/n:\ k\geq 0\}\) and can be simulated exactly
from the recursion
\[
A^{(n)}\bigl((k+1)/n\bigr)=e^{-\kappa/n}A^{(n)}(k/n)+\bigl(1-e^{-\kappa/n}\bigr)X^{(n)}(k/n).
\]
For \(d=2\), \(\kappa=r_{0}=\Tcal=1\), \(\lambda=1\), and increments whose two
coordinates are independent and distributed as \(\mathcal{X}-1\), with
\(\mathcal{X}\) standard exponential, a Monte Carlo experiment
(\(2\times 10^{5}\) paths, standard errors about \(5\times 10^{-4}\)) gives
\(\E[e^{-\tau^{(n)}}]\approx 0.508\), \(0.532\), \(0.549\), \(0.559\) for
\(n=25\), \(100\), \(400\), \(1600\), against the limit
\(\hyp(1/2;1;1)^{-1}\approx 0.570\) in \eqref{eq:kummer-limit}. In this
experiment, the observed errors are consistent with a decay of order
\(n^{-1/2}\); this is an empirical observation only, and we make no general
claim about the rate of convergence.
\end{enumerate}
\end{remark}

\begin{remark}[Domains coupled to the state through a time-discontinuous functional: drawdowns]\label{rem:drawdown}
The augmentation device of Step~2 of the proof of
Proposition~\ref{prop:coupled} is also useful when the natural barrier
of the approximating domains is discontinuous in time. A classical example
is the trailing-stop, or drawdown, domain; the notation of this remark is
local. For the one-dimensional walk \eqref{eq:donsker-walk} and \(h>0\), let
\(\overline{X}^{(n)}(t):=\sup_{s\leq t}X^{(n)}(s)\) and
\(\Gamma^{(n)}_{t}:=(\overline{X}^{(n)}(t)-h,\infty)\), so that the exit time
\(\tau^{(n)}:=\tau(X^{(n)},\Gamma^{(n)})=\inf\{t\geq 0:\ \overline{X}^{(n)}(t)-X^{(n)}(t)\geq h\}\)
is the first time the drawdown of \(X^{(n)}\) reaches \(h\). Since
\(\overline{X}^{(n)}\) is nondecreasing and right-continuous, hence upper
semicontinuous, the admissible space--time set
\(\cO(\Gamma^{(n)})=\{(t,z)\in E:\ z>\overline{X}^{(n)}(t)-h\}\) is open, so that
continuous barriers exist by Lemma~\ref{lem:universality}. However, the
natural barrier \((t,z)\mapsto\overline{X}^{(n)}(t)-h-z\) is discontinuous in
\(t\), and an abstract continuous barrier would make both the convergence
and the verification of \NT{} less transparent. Augmentation gives a direct
route: with \(Z^{(n)}:=(X^{(n)},\overline{X}^{(n)})\in\mathbb{D}^{2}\) and the
deterministic, time-independent barrier \(\Psi_{h}(t,z_{1},z_{2}):=z_{2}-z_{1}-h\),
one has \(\tau^{(n)}=\tau(Z^{(n)},\Psi_{h})\). The map
\(x\mapsto(x,\sup_{s\leq\cdot}x(s))\) is continuous from \(\Done\) to
\(\mathbb{D}^{2}\), because the same time changes in
\eqref{eq:J1-char} serve both components; hence
\(Z^{(n)}\Rightarrow Z:=(B,\overline{B})\), with
\(\overline{B}_{t}:=\sup_{s\leq t}B_{s}\), by Donsker's theorem and the
continuous mapping theorem. The scalarized limit process
\(Y_{t}=\overline{B}_{t}-B_{t}-h\) is continuous, so that \eqref{eq:NTminus}
holds, and \(\tau:=\tau(Z,\Psi_{h})\leq\inf\{t\geq 0:\ B_{t}\leq-h\}<\infty\)
almost surely. At time \(\tau\),
\(B_{\tau}=\overline{B}_{\tau}-h<\overline{B}_{\tau}\), so that \(\overline{B}\)
is constant on \([\tau,\tau+\delta]\) for some random \(\delta>0\). Therefore
\(Y=Y_{0}+M+V\), with \(Y_{0}=-h\), \(M:=-B\) and \(V:=\overline{B}\), satisfies
the assumptions of Proposition~\ref{prop:semimartingale}, \(V\) being
constant on \([\tau,\tau+\delta]\); this gives \eqref{eq:NTplus}, and hence
\NT{}, almost surely. Theorem~\ref{thm:weak} yields
\(\tau^{(n)}\Rightarrow\tau\). By Lévy's theorem~\cite[Section~VI.2]{RevuzYor1999},
\(\overline{B}-B\) is a reflected Brownian motion, so that
\(\E[e^{-\lambda\tau}]=1/{\cosh(h\sqrt{2\lambda})}\) for \(\lambda>0\);
see~\cite{Taylor1975} for drawdowns of Brownian motion with drift
and~\cite{Lehoczky1977} for diffusions. The same argument shows that, for
each fixed level \(u>-h\), \((Z,\Psi_{h})\) satisfies \NT{} at level \(u\)
almost surely. Let \(Y^{(n)}:=\overline{X}^{(n)}-X^{(n)}-h\), so that
\(T^{+}_{Y^{(n)}}(u)\) is the first time at which the drawdown of \(X^{(n)}\)
exceeds \(h+u\). Applying the fixed-level statement to two levels
\(u_{0}<u_{1}\) with \(u_{0}>-h\), we obtain via Theorem~\ref{thm:profiles}
that the profiles \(u\mapsto T^{+}_{Y^{(n)}}(u)\) converge in distribution in
\(\Ms[u_{0},u_{1}]\) to \(u\mapsto T^{+}_{Y}(u)\). Since
\(T^{+}_{Y}(u_{1})\leq\inf\{t\geq 0:\ B_{t}<-h-u_{1}\}<\infty\) almost surely,
the convergence also holds in \((D[u_{0},u_{1}],M_{1})\) in the sense of
Remark~\ref{rem:finite-profiles}.
\end{remark}

\appendix
\renewcommand{\thelemma}{\Alph{section}.\arabic{lemma}}

\section{Proof of Lemma~\ref{lem:monotone-M1}}\label{app:M1}

\subsection{The \texorpdfstring{\(M_{1}\)}{M1} metric via parametric representations}\label{app:M1-metric}

Let \(a<b\). As in Section~\ref{sec:profile-setup}, \(D[a,b]\) denotes the
space of functions \(g:[a,b]\to\R\) that are right-continuous on \([a,b)\)
and have left limits on \((a,b]\); we set \(g(a-):=g(a)\). We only need the
\(M_{1}\) topology on nondecreasing functions, and we recall its definition
in that case; see~\cite{Whitt2002} for the general case. For a
nondecreasing \(g\in D[a,b]\), the completed graph of \(g\) is
\[
G_{g}:=\bigl\{(t,z)\in[a,b]\times\R:\ g(t-)\leq z\leq g(t)\bigr\},
\]
ordered by \((t_{1},z_{1})\preceq(t_{2},z_{2})\) if either \(t_{1}<t_{2}\), or
\(t_{1}=t_{2}\) and \(z_{1}\leq z_{2}\). A parametric representation of \(g\)
is a continuous map \((r,v):[0,1]\to G_{g}\) that is onto and nondecreasing
for \(\preceq\); let \(\Pi(g)\) denote the set of parametric representations
of \(g\). The \(M_{1}\) metric is
\begin{equation}\label{eq:M1-metric}
d_{M_{1}}(g,h):=\inf_{\substack{(r,v)\in\Pi(g)\\(r',v')\in\Pi(h)}}\max\bigl\{\|r-r'\|_{\infty},\|v-v'\|_{\infty}\bigr\},
\end{equation}
where \(\|\cdot\|_{\infty}\) denotes the supremum norm. Since every parametric
representation satisfies \((r(0),v(0))=(a,g(a))\) and
\((r(1),v(1))=(b,g(b))\), \(M_{1}\) convergence implies convergence at the end
points \(a\) and \(b\). By \eqref{eq:M1-metric}, to prove
Lemma~\ref{lem:monotone-M1} it suffices to construct
\((\hat{r},\hat{v})\in\Pi(f)\) and \((\hat{r}_{n},\hat{v}_{n})\in\Pi(f_{n})\)
such that
\begin{equation}\label{eq:M1-suffices}
\|\hat{r}_{n}-\hat{r}\|_{\infty}\to 0\qquad\text{and}\qquad\|\hat{v}_{n}-\hat{v}\|_{\infty}\to 0 .
\end{equation}

\subsection{A canonical parametric representation of a nondecreasing function}\label{app:canonical}

Let \(g\in D[a,b]\) be nondecreasing, and define
\[
L_{g}:=(b-a)+g(b)-g(a)>0,\quad
\Theta_{g}(t):=(t-a)+g(t)-g(a),\quad t\in[a,b].
\]
The function \(\Theta_{g}\) is right-continuous on \([a,b)\), strictly
increasing, with \(\Theta_{g}(a)=0\), \(\Theta_{g}(b)=L_{g}\), and
\(\Theta_{g}(t)-\Theta_{g}(s)\geq t-s\) for \(s\leq t\). For
\(\ell\in[0,L_{g}]\), let
\begin{equation}\label{eq:canonical}
r_{g}(\ell):=\inf\{t\in[a,b]:\ \Theta_{g}(t)\geq\ell\},\quad
v_{g}(\ell):=g(a)+\ell-\bigl(r_{g}(\ell)-a\bigr),
\end{equation}
and let \(\hat{r}_{g}(s):=r_{g}(sL_{g})\) and \(\hat{v}_{g}(s):=v_{g}(sL_{g})\) for
\(s\in[0,1]\). The set in \eqref{eq:canonical} contains \(b\), so \(r_{g}\) is
well defined.

\begin{lemma}\label{lem:canonical}
Let \(g\in D[a,b]\) be nondecreasing. Then \(r_{g}\) is nondecreasing and
\(1\)-Lipschitz on \([0,L_{g}]\), with \(r_{g}(0)=a\) and \(r_{g}(L_{g})=b\), and
\((\hat{r}_{g},\hat{v}_{g})\in\Pi(g)\).
\end{lemma}

\begin{proof}
Clearly \(r_{g}\) is nondecreasing, \(r_{g}(0)=a\), and \(r_{g}(L_{g})=b\)
because \(\Theta_{g}(t)<L_{g}\) for \(t<b\). Let \(0\leq\ell<\ell'\leq L_{g}\),
\(t:=r_{g}(\ell)\) and \(t':=r_{g}(\ell')\), and suppose that \(t'>t\). If
\(t<b\), right-continuity gives \(\Theta_{g}(t)\geq\ell\), and for every
\(s\in(t,t')\) we have \(\Theta_{g}(s)<\ell'\) by the definition of
\(t'\); hence \(s-t\leq\Theta_{g}(s)-\Theta_{g}(t)<\ell'-\ell\), and letting
\(s\uparrow t'\) yields \(t'-t\leq\ell'-\ell\). (The case \(t=b\) is excluded by
\(t'>t\).) Thus \(r_{g}\) is \(1\)-Lipschitz, and consequently \(v_{g}\) is
continuous and nondecreasing, since
\(v_{g}(\ell')-v_{g}(\ell)=(\ell'-\ell)-(r_{g}(\ell')-r_{g}(\ell))\geq 0\).

Let \(\ell\in[0,L_{g}]\) and \(t:=r_{g}(\ell)\). For \(s<t\) we have
\(\Theta_{g}(s)<\ell\), so that \(\Theta_{g}(t-)\leq\ell\) (with
\(\Theta_{g}(a-):=0\)), while \(\Theta_{g}(t)\geq\ell\), by right-continuity if
\(t<b\) and because \(\Theta_{g}(b)=L_{g}\) if \(t=b\). Subtracting \(t-a\) and
adding \(g(a)\) gives \(g(t-)\leq v_{g}(\ell)\leq g(t)\), that is,
\((r_{g}(\ell),v_{g}(\ell))\in G_{g}\). Since both coordinates are
nondecreasing in \(\ell\), the map \(\ell\mapsto(r_{g}(\ell),v_{g}(\ell))\) is
nondecreasing for \(\preceq\).

Finally, the map is onto \(G_{g}\). Let \((t,z)\in G_{g}\) and
\(\ell:=(t-a)+z-g(a)\), so that \(\Theta_{g}(t-)\leq\ell\leq\Theta_{g}(t)\) and
\(\ell\in[0,L_{g}]\). Then \(r_{g}(\ell)\leq t\), and for \(s<t\),
\(\Theta_{g}(s)\leq(s-a)+g(t-)-g(a)<\Theta_{g}(t-)\leq\ell\), so that
\(r_{g}(\ell)=t\) and \(v_{g}(\ell)=z\). The rescaled pair
\((\hat{r}_{g},\hat{v}_{g})\) inherits all these properties, hence belongs to
\(\Pi(g)\).
\end{proof}

\subsection{Convergence of the canonical representations}\label{app:convergence}

\begin{lemma}\label{lem:dini}
Let \(L>0\), let \(h_{n}:[0,L]\to\R\) be nondecreasing, and let
\(h:[0,L]\to\R\) be continuous. If \(h_{n}(\ell)\to h(\ell)\) for every \(\ell\)
in a dense subset of \([0,L]\) containing \(0\) and \(L\), then
\(\|h_{n}-h\|_{\infty}\to 0\).
\end{lemma}

\begin{proof}
Let \(\eta>0\). By uniform continuity, there is \(\delta>0\) such that
\(|h(\ell)-h(\ell')|<\eta/3\) whenever \(|\ell-\ell'|<\delta\). Choose points
\(0=\ell_{0}<\ell_{1}<\cdots<\ell_{m}=L\) of the dense set with
\(\max_{i}(\ell_{i}-\ell_{i-1})<\delta\), and let \(n_{0}\) be such that
\(|h_{n}(\ell_{i})-h(\ell_{i})|<\eta/3\) for all \(i\) and all \(n\geq n_{0}\).
For \(n\geq n_{0}\) and \(\ell\in[\ell_{i-1},\ell_{i}]\), monotonicity gives
\[
h(\ell)-\tfrac{2\eta}{3}<h(\ell_{i-1})-\tfrac{\eta}{3}<h_{n}(\ell_{i-1})\leq h_{n}(\ell)\leq h_{n}(\ell_{i})<h(\ell_{i})+\tfrac{\eta}{3}<h(\ell)+\tfrac{2\eta}{3},
\]
so that \(\|h_{n}-h\|_{\infty}<\eta\).
\end{proof}

\begin{proof}[Proof of Lemma~\ref{lem:monotone-M1}]
Let \(f_{n},f\in D[a,b]\) be nondecreasing, and let \(D_{0}\subset[a,b]\) be
a dense set containing \(a\) and \(b\) such that \(f_{n}(t)\to f(t)\) for all
\(t\in D_{0}\). Write \(L_{n}:=L_{f_{n}}\), \(L:=L_{f}\), \(\Theta_{n}:=\Theta_{f_{n}}\),
\(\Theta:=\Theta_{f}\), \(r_{n}:=r_{f_{n}}\) and \(r:=r_{f}\), and extend \(r_{n}\) to
\([0,\infty)\) by \(r_{n}(\ell):=b\) for \(\ell>L_{n}\); the extended function is
still nondecreasing and \(1\)-Lipschitz, since \(r_{n}(L_{n})=b\).

\emph{Step 1: convergence of the lengths.} Since \(a,b\in D_{0}\),
\(L_{n}=(b-a)+f_{n}(b)-f_{n}(a)\to L\).

\emph{Step 2: pointwise convergence of the inverses.} Clearly
\(r_{n}(0)=a=r(0)\). Let \(\ell\in(0,L)\) and \(t:=r(\ell)\); then \(t>a\),
because \(\Theta(a)=0<\ell\). Let \(\varepsilon>0\). Since \(D_{0}\) is dense,
there is \(t_{-}\in D_{0}\cap(\max\{a,t-\varepsilon\},t)\), and
\(\Theta(t_{-})<\ell\) by the definition of \(r(\ell)\). Since \(t_{-},a\in D_{0}\),
\(\Theta_{n}(t_{-})\to\Theta(t_{-})\), so that \(\Theta_{n}(t_{-})<\ell\) for all
large \(n\); as \(\Theta_{n}\) is increasing, \(\Theta_{n}(s)<\ell\) for all
\(s\leq t_{-}\), whence \(r_{n}(\ell)\geq t_{-}>t-\varepsilon\) (this is also true
if \(\ell>L_{n}\)). If \(t=b\), then \(r_{n}(\ell)\leq b=t\). If \(t<b\), there
is \(t_{+}\in D_{0}\cap(t,\min\{b,t+\varepsilon\})\), and
\(\Theta(t_{+})>\Theta(t)\geq\ell\); hence \(\Theta_{n}(t_{+})>\ell\) for all
large \(n\), and \(r_{n}(\ell)\leq t_{+}<t+\varepsilon\). Thus
\(r_{n}(\ell)\to r(\ell)\) for every \(\ell\in[0,L)\). Finally,
\(r_{n}(L)\leq b=r(L)\), and \(r_{n}(L)\geq r_{n}(\ell')\to r(\ell')\) for every
\(\ell'<L\), while \(r(\ell')\to r(L)\) as \(\ell'\uparrow L\) by continuity of
\(r\); hence \(r_{n}(L)\to r(L)\) as well.

\emph{Step 3: uniform convergence.} By Step~2 and Lemma~\ref{lem:dini},
applied with \(h_{n}:=r_{n}|_{[0,L]}\) and \(h:=r\), which is continuous by
Lemma~\ref{lem:canonical}, we obtain \(\sup_{\ell\in[0,L]}|r_{n}(\ell)-r(\ell)|\to 0\).

\emph{Step 4: convergence of the rescaled representations.} We write
\(\hat{r}_{n}:=\hat{r}_{f_{n}}\), \(\hat{v}_{n}:=\hat{v}_{f_{n}}\),
\(\hat{r}:=\hat{r}_{f}\) and \(\hat{v}:=\hat{v}_{f}\). For \(s\in[0,1]\), the
\(1\)-Lipschitz property of the extended \(r_{n}\) gives
\[
\begin{aligned}
|\hat{r}_{n}(s)-\hat{r}(s)|
&\leq|r_{n}(sL_{n})-r_{n}(sL)|+|r_{n}(sL)-r(sL)|\\
&\leq|L_{n}-L|+\sup_{\ell\in[0,L]}|r_{n}(\ell)-r(\ell)|,
\end{aligned}
\]
and, by \eqref{eq:canonical},
\[
|\hat{v}_{n}(s)-\hat{v}(s)|\leq|f_{n}(a)-f(a)|+|L_{n}-L|+|\hat{r}_{n}(s)-\hat{r}(s)| .
\]
By Steps~1 and~3 and since \(a\in D_{0}\), both bounds tend to \(0\) uniformly
in \(s\in[0,1]\). Since \((\hat{r}_{n},\hat{v}_{n})\in\Pi(f_{n})\) and
\((\hat{r},\hat{v})\in\Pi(f)\) by Lemma~\ref{lem:canonical},
\eqref{eq:M1-suffices} holds, and \(d_{M_{1}}(f_{n},f)\to 0\).
\end{proof}

\section*{Acknowledgements}

The author would like to thank the Editors and the referees for their
constructive feedback and their suggestions that helped improve the initial
manuscript substantially.

\begingroup
\raggedright

\endgroup

\end{document}